\documentclass[11pt,a4paper]{article}

\usepackage[T1]{fontenc}
\usepackage[utf8]{inputenc}
\usepackage{lmodern}
\usepackage[margin=1in]{geometry}
\usepackage{indentfirst}
\usepackage{amsmath,amssymb,amsthm,mathtools}
\usepackage{graphicx}
\usepackage{microtype}
\usepackage[numbers,sort&compress]{natbib}
\usepackage[
  colorlinks=true,
  linkcolor=blue,  
  citecolor=blue,  
  urlcolor=blue   
]{hyperref}

\newtheorem{theorem}{Theorem}[section]
\newtheorem{proposition}[theorem]{Proposition}

\theoremstyle{definition}

\newcommand{\dd}{\mathrm{d}}

\title{A Discrete Curvature Action for Origami Dynamical Systems}

\author{Ryutaro Ichikawa, Mitsuru Shibayama}

\date{}

\begin{document}

\maketitle

\begin{abstract}
Nonuniform folded configurations in rigid origami can be described by a dynamical system governing spatial propagation between closed boundaries. In this work, we close the folded surface with oriented auxiliary faces and use the contribution of the creases and boundary edges to the signed total discrete mean curvature of the resulting polyhedral surface to construct a curvature action for two-dimensional dynamical systems of tubular origami. Using the Schl\"afli formula, we show that the differential of this action is expressed solely in terms of the input and output boundaries. The same construction applies to Maxwell origami tubes with unequal left and right edge lengths and no assumed rotational symmetry of the folded configuration, as well as to non-tubular polygonal Spidrons in which propagation proceeds inward with a change of scale.
Under a local twist condition, the sequences of boundary angles that determine the overall folded shape are characterized variationally as stationary points of the discrete action with the endpoint angles held fixed. By imposing geometric constraints, we numerically construct origami surfaces with distinctive shapes, including knotted tubes, opening the way to the exploration of new folded configurations.
\end{abstract}

\section{Introduction}
\label{sec:introduction}

Origami tessellations with periodic crease patterns have been widely studied as structures whose overall shape and deformation can be designed through local geometry. Schenk and Guest related the geometry of Miura-ori to its deformation characteristics and demonstrated its use in mechanical metamaterials~\cite{Schenk2013Geometry}. McInerney et al.\ identified hidden symmetries underlying rigid folding mechanisms in periodic origami~\cite{McInerney2020Hidden}. Nassar et al.\ derived constraints on the metric and curvature of macroscopic surfaces approximated by the eggbox pattern from its unit-scale geometry~\cite{Nassar2017Curvature}. These studies provide a foundation for understanding the geometry and deformation of repetitive folded structures. However, a periodic crease pattern need not produce a folded configuration with the same periodicity. Neighboring units can occupy different states while remaining geometrically compatible, giving rise to spatially nonuniform configurations.

To describe such configurations, Imada and Tachi formulated tubular origami as spatial propagation from one boundary to the next. For a waterbomb tube under rotational symmetry assumptions, edge-length constraints and prescribed mountain--valley assignments determine the state of a neighboring ring from the current ring state. The configuration of the entire tube is then represented by an orbit of a discrete dynamical system~\cite{Imada2022Waterbomb}. They subsequently extended this approach to a broader family of rotationally symmetric tubular origami tessellations, establishing a formulation in terms of two-dimensional area-preserving maps and also treating nonconservative systems involving changes of scale~\cite{Imada2023Undulations}. Here, iteration describes spatial propagation along the structure rather than the temporal evolution of a folding motion.

This dynamical-systems approach extends beyond rotationally symmetric tubular origami. In Maxwell lattices, the numbers of kinematic variables and constraints balance within each unit cell. The relationship between deformations localized near a boundary and those extending into the bulk has been studied in this setting, and related ideas have been introduced into origami and kirigami~\cite{Kane2014Topological,Chen2016Topological}. Building on this perspective, the Maxwell origami tube is described by a higher-dimensional propagation map that accommodates multiple boundary degrees of freedom without imposing rotational symmetry on the folded configuration~\cite{Imada2025Maxwell}. Dynamical-systems descriptions have also been developed for Spidron-inspired origami, in which units at different scales are connected recursively. These systems exhibit different nonlinear behaviors depending on the folding mode~\cite{Imada2026Spidron}. Together, these studies show how origami structures with different symmetries, scales, and shapes can be understood through spatial propagation governed by local geometry.

A central question concerns the geometric structure underlying these propagation maps. Area preservation has been established for rotationally symmetric tubular origami~\cite{Imada2023Undulations}. For the Maxwell origami tube, Imada and Tachi~\cite{Imada2025Maxwell} observed numerically that the eigenvalues of the linearized propagation map at fixed points occur in reciprocal pairs, suggesting an underlying conservative structure. Wherever the boundary relation is locally the graph of a smooth propagation map, Miranda's boundary isotropy theorem~\cite{Miranda2025Domes} implies that the map is symplectic. As discussed in Subsection~\ref{subsec:maxwell-ring-action}, the reciprocal pairing at fixed points follows from this symplecticity.

Symplecticity provides a starting point for a variational description of propagation. Under a suitable twist condition, the map admits a local generating function. Summing this function over consecutive pairs of states defines a discrete action, and the discrete variational principle characterizes propagation orbits locally by stationary sequences with fixed endpoint angles~\cite{Meiss1992}. In billiards, for example, this action is the path length, and stationarity with respect to the reflection points yields the reflection law~\cite{Tabachnikov2005Billiards}. For origami propagation, however, the existence of a local generating function leaves open the nontrivial problem of identifying a quantity of the folded geometry that gives it. Here we identify a quantity that plays the role of path length in spatial origami propagation.

In our previous work~\cite{Ichikawa2026InvariantCurves}, we obtained a curvature interpretation of the generating function in the formal integrable limit of rotationally symmetric tubular origami under edge-length symmetry assumptions. There, the generating function was obtained by explicitly integrating the propagation relations, and the Schl\"afli formula was used to verify the generating relations. The question addressed here is how to construct the curvature action directly from geometry, without first integrating the propagation map.

The key observation is that closing the surface of each propagation step with oriented auxiliary faces turns the Schläfli formula~\cite{Milnor1994Schlafli,Regge1961,RivinSchlenker1999} into a variational relation between boundary variables. The lengths of the physical creases remain fixed, whereas the diagonal lengths of the auxiliary surfaces spanning the boundaries vary with the folded configuration. Separating these contributions shows that the fixed-length-edge contribution to total discrete mean curvature provides an action relating the input and output boundary variables. This allows a generating function to be constructed geometrically, under a local twist condition, without integrating an explicit expression for the propagation map.

Rotational symmetry reduces the number of boundary degrees of freedom, but the boundary variation formula also applies to folded configurations without rotational symmetry. For the Maxwell origami tube extended to unequal prescribed left and right edge lengths, we use the half-length of each cap diagonal and its signed exterior dihedral angle as a conjugate pair of boundary variables. These pairs give a variational description of the higher-dimensional propagation. The same choice of boundary variables applies to inward propagation in Spidron-type origami generalized to polygonal boundaries, extending the variational structure to a non-tubular system. When successive rings differ in scale, the normalized propagation is conformally symplectic, and the scaling of the curvature sum determines the weights in the discrete action. The boundary coordinates and the curvature action together link local folding compatibility to the shape of the entire structure. This relationship also suggests how to search for folded configurations
by starting from their overall shape. For Maxwell origami tubes, for example,
we place the rings along a target knot to form an initial configuration
and seek a folded shape that satisfies the edge-length constraints and
closes in space. For polygonal Spidrons, we can choose the shape of the
outer boundary and examine the folded configurations that arise through
inward propagation. Folded shapes obtained under these different geometric
conditions can thus be treated within the same variational framework
based on the curvature action. This opens a way to explore new folded
shapes beyond tubular origami.

The paper is organized as follows. Section~\ref{sec:map} introduces
tubular origami propagation and its local generating-function formulation.
Section~\ref{sec:action} develops the curvature action, its composition,
and its extension to changes of scale. Section~\ref{sec:maxwell} treats
the generalized Maxwell origami tube, and Section~\ref{sec:spidron-boundary}
considers inward propagation in polygonal Spidrons.
Section~\ref{sec:numerical-examples} presents numerical constructions
of folded configurations for both systems, including knotted Maxwell
origami tubes, and relates these configurations to stationary points
of the corresponding discrete actions. Finally, Section~\ref{sec:conclusion}
presents our conclusions and discusses open problems.

\section{Tubular origami and the propagation map}
    \label{sec:map}

\begin{figure}[tb]
      \centering
      \includegraphics[width=\textwidth]{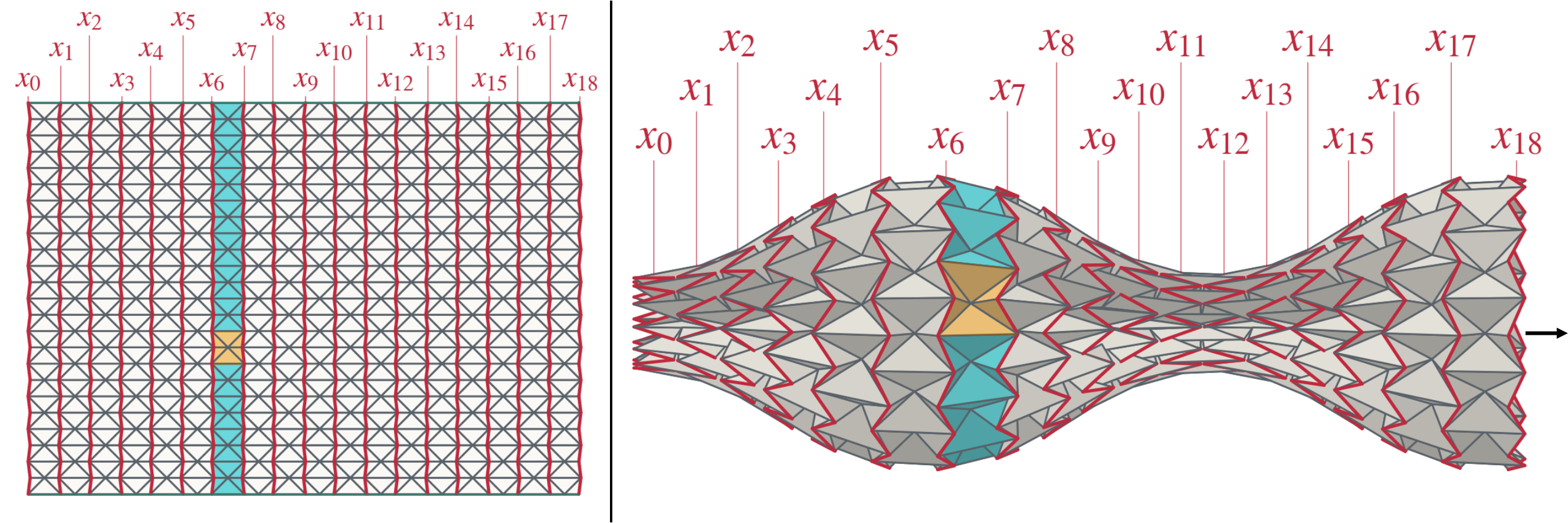}
      \caption{
        Tubular origami viewed as spatial propagation under the map $M$.
        Left: the planar crease pattern, whose upper and lower
        boundaries (green) are identified to form a tube.
        Right: the corresponding folded configuration.
        Red lines indicate corresponding circumferential zigzags
        in successive rings, with
        $\mathbf{x}_n=(\theta_n,I_n)^{\mathsf T}$
        denoting the state vector at the $n$th ring and
        $\mathbf{x}_{n+1}=M(\mathbf{x}_n)$.
        The cyan region indicates a ring, and the orange region
        indicates one module within the ring.
        The arrow indicates the direction of spatial propagation.
      }
      \label{fig:intro-origami}
    \end{figure}

    We adopt the tubular origami model of Imada and Tachi
    \cite{Imada2023Undulations}.
    A tube is formed by identifying corresponding points on the
    upper and lower boundaries of an origami tessellation.
    We refer to a repeating unit of the crease pattern as a
    \emph{module}.
    As illustrated in Fig.~\ref{fig:intro-origami},
    $N\in\mathbb{Z}_{\geq 3}$ copies of the orange module
    are arranged around the circumference to form a \emph{ring},
    highlighted in cyan.
    We consider folded configurations with $N$-fold rotational
    symmetry about the tube axis.
    
    Under this rotational symmetry, each circumferential \emph{zigzag}
    is generated by rotating a representative segment contained
    in a single module through multiples of $2\pi/N$.
    We consider rigid origami, in which crease lengths and
    the shape of each face remain unchanged during folding.
    Under these constraints, the configuration of the representative
    segment, and hence of the entire zigzag, is described by
    the state vector $\mathbf{x}=(\theta,I)^{\mathsf T}$
    \cite{Ichikawa2026InvariantCurves}.
    Geometric compatibility conditions between neighboring zigzags,
    together with the mountain--valley assignments, determine
    the propagation maps between successive zigzags.
    Composing the elementary maps within a representative module
    gives the two-dimensional map $M$ from one ring to the next.
    Writing $\mathbf{x}_n=(\theta_n,I_n)^{\mathsf T}$
    for the state vector of the corresponding zigzag in the $n$th ring,
    the folded configuration of the tube is described by a sequence
    satisfying $\mathbf{x}_{n+1}=M(\mathbf{x}_n)$.
    
    Throughout the paper, signed angular variables are regarded as elements of
    $\mathbb R/(2\pi\mathbb Z)$. When principal representatives are needed,
    we use $[-\pi,\pi)$. In differential and variational calculations,
    we use locally smooth real-valued lifts, denoted by the same symbols.
    Prescribed geometric parameters and the selected folding solution
    are held fixed unless stated otherwise.

    \subsection{Zigzag parameterization and local propagation}
        \label{subsec:zigzag-propagation}
        
\begin{figure}[tb]
          \centering
          \includegraphics[width=\textwidth]{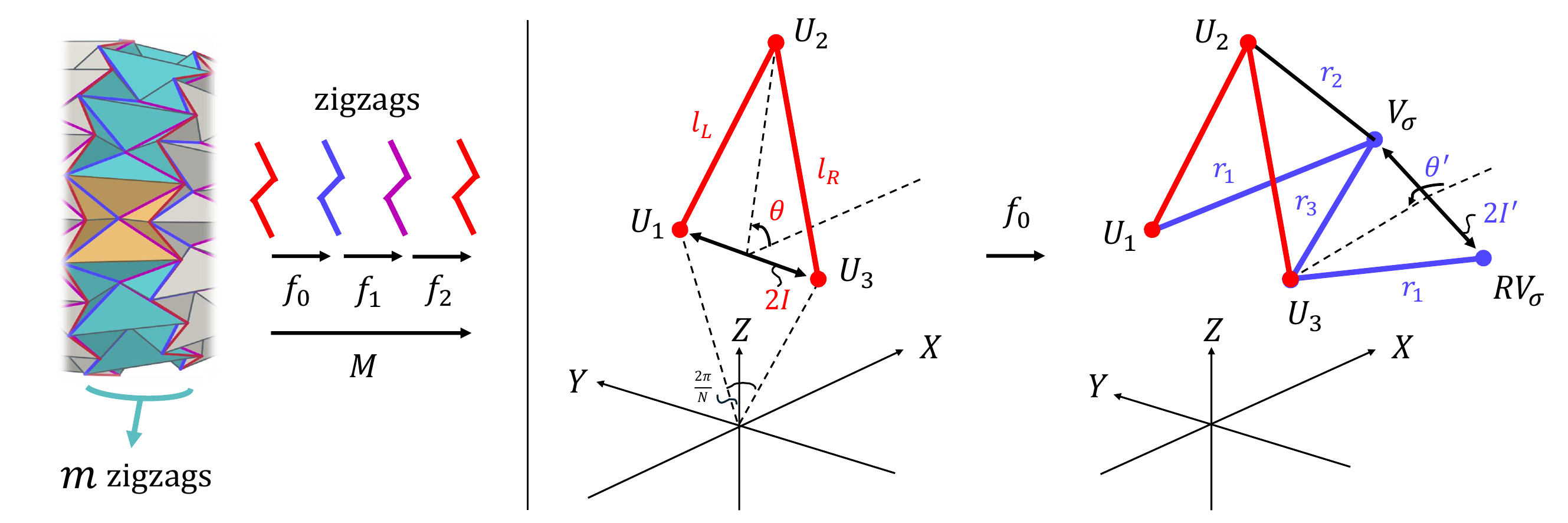}
          \caption{
            Zigzag parameterization and elementary propagation.
            Left: decomposition of a ring into circumferential zigzags
            and propagation between successive zigzags.
            The illustrated example has $m=3$ and requires no
            reparameterization, so $M=f_2\circ f_1\circ f_0$.
            The final red zigzag belongs to the next ring.
            Center: the input segment $(U_1,U_2,U_3)$,
            with crease lengths $l_L,l_R$ and state vector
            $\mathbf{x}=(\theta,I)^{\mathsf T}$.
            Here, $2I=|U_1U_3|$, and $\theta$ is the signed angle
            from the positive $X$-direction to the direction perpendicular
            to the line $U_1U_3$ and pointing toward $U_2$,
            with positive angles taken toward the positive $Z$-direction.
            Right: construction of $V_\sigma$ from the distance constraints
            $|U_jV_\sigma|=r_j$, $j=1,2,3$.
            The output segment is $(V_\sigma,U_3,V_\sigma')$,
            where $V_\sigma'$ is obtained by rotating $V_\sigma$
            about the tube axis through $2\pi/N$.
            Its state vector $\mathbf{x}'=(\theta',I')^{\mathsf T}$
            is defined by the same convention, with
            $2I'=|V_\sigma V_\sigma'|$.
            The $X$-axis is the tube axis.
          }
          \label{fig:zigzag-parametrization}
        \end{figure}

        Each ring contains $m\in\mathbb{Z}_{\geq 1}$ circumferential
        zigzags arranged in the direction of spatial propagation.
        To describe propagation between successive zigzags,
        we introduce the elementary map
        \begin{equation}
          \mathbf{x}'
          =
          f(\mathbf{x};\,N,l_L,l_R,r_1,r_2,r_3,\sigma),
          \qquad
          \mathbf{x}=(\theta,I)^{\mathsf T},
          \quad
          \mathbf{x}'=(\theta',I')^{\mathsf T}.
          \label{eq:elementary-propagation}
        \end{equation}
        Here, $l_L,l_R,r_1,r_2,r_3\in\mathbb{R}_{>0}$
        are prescribed crease lengths, and $\sigma\in\{+1,-1\}$
        selects the sphere-intersection point consistent with the prescribed
        mountain--valley assignment in the geometric construction described below.
        Figure~\ref{fig:zigzag-parametrization} illustrates
        the geometric construction of $f$ and the composition
        of elementary steps over one module.
        
        Let $U_1,U_2,U_3$ be three consecutive vertices of
        the input zigzag, with
        \begin{equation*}
          |U_1U_2|=l_L,
          \qquad
          |U_2U_3|=l_R.
        \end{equation*}
        Choose the tube axis as the $X$-axis, with its positive
        direction chosen so that the right-handed rotation
        through $2\pi/N$ maps $U_1$ to $U_3$.
        The origin is the intersection of the tube axis with
        the plane perpendicular to it containing $U_1$ and $U_3$.
        Let $\mathbf{u}_j\in\mathbb{R}^3$ denote the position vector
        of $U_j$, for $j=1,2,3$.
        Then
        \begin{equation*}
          \mathbf{u}_3=R\mathbf{u}_1,
          \qquad
          R=R_X(2\pi/N),
        \end{equation*}
        where $R_X(\alpha)$ denotes the right-handed rotation matrix
        through angle $\alpha$ about the $X$-axis.
        We choose the positive $Y$-direction from $U_3$ toward $U_1$
        and complete a right-handed orthonormal frame by setting
        \begin{equation*}
          \mathbf{e}_Y
          =
          \frac{\mathbf{u}_1-\mathbf{u}_3}
               {\|\mathbf{u}_1-\mathbf{u}_3\|},
          \qquad
          \mathbf{e}_Z=\mathbf{e}_X\times\mathbf{e}_Y,
        \end{equation*}
        where $\mathbf{e}_X$ is the unit vector along the
        positive $X$-axis.
        
        Define
        \begin{equation*}
          I=\frac12|U_1U_3|.
        \end{equation*}
        In this coordinate frame, the position vectors of the
        endpoints are
        \begin{equation*}
          \mathbf{u}_1
          =
          \bigl(0,\,I,\,I\cot(\pi/N)\bigr)^{\mathsf T},
          \qquad
          \mathbf{u}_3
          =
          \bigl(0,\,-I,\,I\cot(\pi/N)\bigr)^{\mathsf T}.
        \end{equation*}
        Let $U_M$ be the orthogonal projection of $U_2$ onto
        the line $U_1U_3$, and let $\mathbf{u}_M$ be its
        position vector.
        The angular variable $\theta\in\mathbb R/(2\pi\mathbb Z)$
        is defined by
        \begin{equation*}
          \frac{\mathbf{u}_2-\mathbf{u}_M}
               {\|\mathbf{u}_2-\mathbf{u}_M\|}
          =
          \cos\theta\,\mathbf{e}_X+\sin\theta\,\mathbf{e}_Z.
        \end{equation*}
        Thus, $\theta$ measures the signed angle from the positive
        $X$-direction toward the positive $Z$-direction in the
        plane perpendicular to $U_1U_3$.
        The input triangle is nondegenerate precisely when
        \begin{equation}
          \frac{|l_L-l_R|}{2}
          <
          I
          <
          \frac{l_L+l_R}{2}.
          \label{eq:input-domain}
        \end{equation}
        Together with the prescribed lengths $l_L,l_R$,
        the state vector $\mathbf{x}=(\theta,I)^{\mathsf T}$
        determines the input segment.
        
        The next vertex $V_\sigma$ is constructed from
        \begin{equation}
          |U_jV_\sigma|=r_j,
          \qquad j=1,2,3.
          \label{eq:sphere-intersection}
        \end{equation}
        These equations describe the intersection of three spheres
        centered at $U_j$ with radii $r_j$.
        When two intersection points exist,
        $\sigma\in\{+1,-1\}$ selects the point corresponding to
        the prescribed mountain or valley fold along $U_2V_\sigma$.
        We work locally where the selected sphere-intersection point
        depends smoothly on the input state and the input and output
        zigzag triangles are nondegenerate.
        
        Let $\mathbf{v}_\sigma$ be the position vector of $V_\sigma$,
        and let $V_\sigma'$ be the vertex with position vector
        $\mathbf{v}_\sigma'=R\mathbf{v}_\sigma$.
        The map $f$ assigns the output state vector
        $\mathbf{x}'=(\theta',I')^{\mathsf T}$
        to the zigzag segment $(V_\sigma,U_3,V_\sigma')$.
        These variables are defined by the same construction as
        $(\theta,I)$, using a local coordinate frame adapted to
        the output segment while retaining the positive axial direction.
        In particular,
        \begin{equation*}
          I'=\frac12|V_\sigma V_\sigma'|.
        \end{equation*}
        The two crease lengths of the output segment are
        \begin{equation*}
          |V_\sigma U_3|=r_3,
          \qquad
          |U_3V_\sigma'|=r_1,
        \end{equation*}
        where the second equality follows from rotational symmetry.
        For connectivity with a circumferential offset,
        this segment is used directly as the input to the next
        propagation step.
        
        Depending on the connectivity pattern, however, the next
        step instead uses the segment of the same zigzag obtained
        by rotating $(U_1,V_\sigma,U_3)$ about the tube axis
        through $2\pi/N$ \cite{Imada2023Undulations}.
        By rotational symmetry, we may equivalently use
        $(U_1,V_\sigma,U_3)$ itself.
        We denote the corresponding reparameterization by
        \begin{equation*}
          \mathbf{x}^*
          =
          g(\mathbf{x}';\,N,r_3,r_1),
          \qquad
          \mathbf{x}^*=(\theta^*,I^*)^{\mathsf T}.
        \end{equation*}
        This operation leaves the folded zigzag unchanged but
        changes the segment used to define its state variables,
        interchanging the left and right crease lengths.
        The state vector $\mathbf{x}^*$ is then used as the input
        to the next propagation step.
        
        To include both connectivity patterns, let
        $k_i\in\{0,1\}$ indicate whether reparameterization is
        required at the $i$th step, where $i=0,\ldots,m-1$.
        We write $f_i$ for the elementary map with the crease
        lengths and mountain--valley assignment of that step,
        and $g_i$ for the corresponding reparameterization.
        The step map is
        \begin{equation*}
          M_i=
          \begin{cases}
            f_i, & k_i=0,\\
            g_i\circ f_i, & k_i=1.
          \end{cases}
        \end{equation*}
        The final step $M_{m-1}$ connects the last zigzag of
        the current ring to the first zigzag of the next ring.
        By periodicity of the crease pattern, the ordered pair of
        crease lengths at the output agrees with that at the input.
        The tubular origami map is therefore
        \begin{equation}
          M=M_{m-1}\circ\cdots\circ M_0.
          \label{eq:module-map}
        \end{equation}
        One iteration of $M$ advances the boundary state from ring
        $n$ to ring $n+1$ according to
        \begin{equation}
          \mathbf{x}_{n+1}=M(\mathbf{x}_n),
          \qquad
          \mathbf{x}_n=(\theta_n,I_n)^{\mathsf T},
          \label{eq:origami-evolution}
        \end{equation}
        where the states are taken at corresponding zigzags.
        When no reparameterization is required,
        $M=f_{m-1}\circ\cdots\circ f_0$.

        \subsection{Generating functions and the discrete variational principle}
        \label{subsec:variational-principle}
        
        The elementary propagation map $f$ and the reparameterization
        map $g$ both preserve the symplectic form
        $\omega=\dd\theta\wedge\dd I$ \cite{Imada2023Undulations}.
        Their compositions $M_i$, and hence the tubular origami map $M$,
        therefore preserve the same form.
        For $\mathbf{x}_{n+1}=M(\mathbf{x}_n)$, this gives
        \begin{equation}
          \dd\theta_{n+1}\wedge\dd I_{n+1}
          =
          \dd\theta_n\wedge\dd I_n.
          \label{eq:area-preservation}
        \end{equation}
        We work locally with smooth real-valued representatives
        of the angular variables.
        
        Regard $\theta_{n+1}$ and $I_{n+1}$ as functions of
        the input variables $(\theta_n,I_n)$, and define the one-form
        $\eta=I_n\,\dd\theta_n-I_{n+1}\,\dd\theta_{n+1}$.
        By Eq.~\eqref{eq:area-preservation}, $\dd\eta=0$.
        The Poincar\'e lemma therefore guarantees, on a sufficiently
        small contractible neighborhood of the input state,
        a smooth function $S(\theta_n,I_n)$ satisfying $\dd S=\eta$.
        
        To express this primitive in terms of consecutive angular
        variables, assume that $M$ satisfies the twist condition
        $\partial\theta_{n+1}/\partial I_n\neq0$.
        By the inverse function theorem, $I_n$ can then be expressed
        locally as a smooth function $I_n=I_n(\theta_n,\theta_{n+1})$.
        Defining
        $L(\theta_n,\theta_{n+1})
        =S\bigl(\theta_n,I_n(\theta_n,\theta_{n+1})\bigr)$,
        we obtain
        \begin{equation*}
          \dd L
          =
          I_n\,\dd\theta_n-I_{n+1}\,\dd\theta_{n+1},
        \end{equation*}
        where both $I_n$ and $I_{n+1}$ are now regarded as
        functions of $(\theta_n,\theta_{n+1})$.
        Let $D_1$ and $D_2$ denote differentiation with respect to the
        first and second arguments of $L$, respectively. For a vector
        argument, $D_iL$ denotes the corresponding column gradient.
        Comparing coefficients yields
        \begin{equation}
          I_n
          =
          D_1L(\theta_n,\theta_{n+1}),
          \qquad
          I_{n+1}
          =
          -D_2L(\theta_n,\theta_{n+1}).
          \label{eq:generating-relations}
        \end{equation}
        Hence, $L$ is a local generating function for $M$.
        The primitive $S$ exists without the twist condition;
        the twist condition is needed to use
        $(\theta_n,\theta_{n+1})$ as local coordinates.
        
        For an angular sequence $\theta_0,\ldots,\theta_J$,
        with $J\geq1$, whose consecutive pairs lie in the domain of $L$,
        define the discrete action
        \begin{equation}
          \mathcal{A}(\theta_0,\ldots,\theta_J)
          =
          \sum_{n=0}^{J-1}L(\theta_n,\theta_{n+1}).
          \label{eq:discrete-action}
        \end{equation}
        For variations with fixed endpoints,
        $\delta\theta_0=\delta\theta_J=0$,
        the stationarity condition $\delta\mathcal{A}=0$
        is equivalent to the discrete Euler--Lagrange equations
        \begin{equation}
          D_2L(\theta_{n-1},\theta_n)
          +
          D_1L(\theta_n,\theta_{n+1})
          =0,
          \qquad n=1,\ldots,J-1.
          \label{eq:discrete-euler-lagrange}
        \end{equation}
        By Eq.~\eqref{eq:generating-relations}, these equations
        require the value of $I_n$ obtained from the preceding step
        to agree with the value used in the following step.
        Thus, after reconstructing $I_n$ through
        Eq.~\eqref{eq:generating-relations}, stationary angular
        sequences correspond locally to orbit segments satisfying
        $\mathbf{x}_{n+1}=M(\mathbf{x}_n)$.
        This is the discrete variational principle for $M$
        \cite{Meiss1992}.
        
        The argument above establishes the local existence of a
        generating function but does not identify its geometric meaning.
        In the next section, we construct such a function for the
        elementary propagation map $f$ from the fixed-length-edge contribution
        to the total discrete mean curvature
        of an associated polyhedral surface.

    \section{Geometric construction of the discrete curvature action}
    \label{sec:action}
    
    \subsection{The Schl\"afli formula}
    \label{subsec:schlafli}

    In a three-dimensional space of constant sectional curvature $\kappa$,
    the classical Schl\"afli formula relates the variation of the volume $V$
    enclosed by a polyhedron to the variations of its interior dihedral
    angles $\alpha_e$~\cite{Milnor1994Schlafli,Regge1961,RivinSchlenker1999}:
    \begin{equation*}
      \sum_e \ell_e\,\dd\alpha_e = 2\kappa\,\dd V,
    \end{equation*}
    where $\ell_e$ is the length of edge $e$. In Euclidean space, the
    right-hand side vanishes.

    Souam~\cite{Souam2004} proves a Schl\"afli formula for oriented polyhedra
    that may self-intersect. In the Euclidean case, his formulation excludes
    antiparallel oriented normals along a shared edge. Such a pair can occur
    where a cap meets an auxiliary face in our construction. The signed
    exterior dihedral angle at such an edge is $\pi$ modulo $2\pi$, so one
    can still choose a smooth real-valued lift locally through that
    configuration. The following proof gives the Euclidean formula with
    this angle convention.

    \begin{proposition}[Euclidean Schl\"afli formula]
    \label{prop:euclidean-schlafli}
    Let $P(t)$ be a smooth family of realizations in $\mathbb R^3$ of a
    finite, closed, oriented polyhedral surface with fixed combinatorial
    structure. Assume that each edge has positive length and each face
    remains planar and nondegenerate. Each edge has two incident faces,
    which induce opposite directions along it. The realization may
    self-intersect, and incidences are determined by the prescribed gluing.

    For an edge $e$ with incident faces $F_1,F_2$, let
    $\boldsymbol{\nu}_1,\boldsymbol{\nu}_2$ be their oriented unit normals
    and let $\mathbf t_e$ be the unit tangent along $e$ in the direction
    induced by $F_1$. Define the signed exterior dihedral angle modulo
    $2\pi$ by
    \begin{equation}
      \cos\Phi_e=\boldsymbol{\nu}_1\cdot\boldsymbol{\nu}_2,
      \qquad
      \sin\Phi_e=\mathbf t_e\cdot
        (\boldsymbol{\nu}_1\times\boldsymbol{\nu}_2).
      \label{eq:signed-dihedral-angle}
    \end{equation}
    For any locally smooth real-valued lifts of these angles,
    \begin{equation*}
      \sum_{e\in E(P)}\ell_e(t)\frac{\dd\Phi_e(t)}{\dd t}=0.
    \end{equation*}
    \end{proposition}

    \begin{proof}
    For each edge, the cosine and sine in the angle definition give a
    smooth map into the unit circle and hence admit a smooth local lift.
    Set $\mathbf b_e=\mathbf t_e\times\boldsymbol{\nu}_1$. In the moving
    frame $(\boldsymbol{\nu}_1,\mathbf b_e)$ perpendicular to $\mathbf t_e$,
    \[
      \boldsymbol{\nu}_2
      =\cos\Phi_e\,\boldsymbol{\nu}_1+\sin\Phi_e\,\mathbf b_e.
    \]
    Differentiating with respect to $t$ and taking the component along
    $\mathbf t_e\times\boldsymbol{\nu}_2$ gives
    \[
      \dot{\boldsymbol{\nu}}_2\cdot(\mathbf t_e\times\boldsymbol{\nu}_2)
      =\dot\Phi_e+\dot{\boldsymbol{\nu}}_1\cdot\mathbf b_e.
    \]
    Thus, with $\boldsymbol{\Omega}_F
    =\boldsymbol{\nu}_F\times\dot{\boldsymbol{\nu}}_F$ for each face $F$,
    \[
      \dot\Phi_e=\mathbf t_e\cdot
      (\boldsymbol{\Omega}_{F_2}-\boldsymbol{\Omega}_{F_1}).
    \]
    This computation is valid at both $\Phi_e\equiv 0$ and
    $\Phi_e\equiv\pi$ modulo $2\pi$.

    Let $\mathbf t_{F,e}$ be the unit tangent in the direction induced
    by $F$ along $e$. Collecting the edge contributions by face, using
    the opposite directions induced along each shared edge, yields
    \[
      \sum_e\ell_e\dot\Phi_e
      =-\sum_F\boldsymbol{\Omega}_F\cdot
        \left(\sum_{e\subset\partial F}\ell_e\mathbf t_{F,e}\right)=0,
    \]
    since the oriented edge vectors around each face sum to zero.
    \end{proof}

    Exchanging $F_1$ and $F_2$ reverses both $\mathbf t_e$ and the cross
    product in the sine formula, leaving the angle unchanged. For a
    convex polyhedron with outward-oriented normals, this angle agrees
    modulo $2\pi$ with the usual exterior angle $\pi-\alpha_e$. On a
    connected neighborhood, two smooth lifts differ by a constant integer
    multiple of $2\pi$ and therefore have the same differential.

    Applying Proposition~\ref{prop:euclidean-schlafli} along each parameter
    direction gives the differential identity
    \begin{equation}
      \sum_{e\in E(P)}\ell_e\,\dd\Phi_e=0.
      \label{eq:schlafli-differential}
    \end{equation}
    In the next subsection, we construct the closed oriented polyhedral
    surface associated with the elementary propagation map $f$.
    Separating the fixed-length edge terms from the terms associated with
    the input and output states then yields the discrete curvature action.

            \subsection{Construction of the oriented polyhedral surface}
            \label{subsec:polyhedron}
            
            We construct a closed oriented polyhedral surface associated
            with one application of the elementary propagation map
            $\mathbf{x}'=f(\mathbf{x})$.
            We use the vertices and state variables introduced in
            Section~\ref{subsec:zigzag-propagation}, with
            $R=R_X(2\pi/N)$, $RU_1=U_3$, and $RV_\sigma=V_\sigma'$.
            All triangular faces introduced below are assumed
            nondegenerate, with $I,I'>0$.
            
            Let $[A,B,C]$ denote the triangular face oriented by
            the indicated vertex ordering.
            For position vectors $\mathbf{a},\mathbf{b},\mathbf{c}$,
            its oriented unit normal is
            \begin{equation*}
              \boldsymbol{\nu}(A,B,C)
              =
              \frac{
                (\mathbf{b}-\mathbf{a})\times(\mathbf{c}-\mathbf{a})
              }{
                \|(\mathbf{b}-\mathbf{a})\times(\mathbf{c}-\mathbf{a})\|
              }.
            \end{equation*}
            Writing $[A,B]$ for the oriented edge from $A$ to $B$,
            the boundary is
            $\partial[A,B,C]=[A,B]+[B,C]+[C,A]$.
            We use the same bracket notation for convex planar
            polygonal faces, oriented by the cyclic order
            of their boundary vertices.
            
            \begin{figure}[tb]
              \centering
              \includegraphics[width=\textwidth]
                {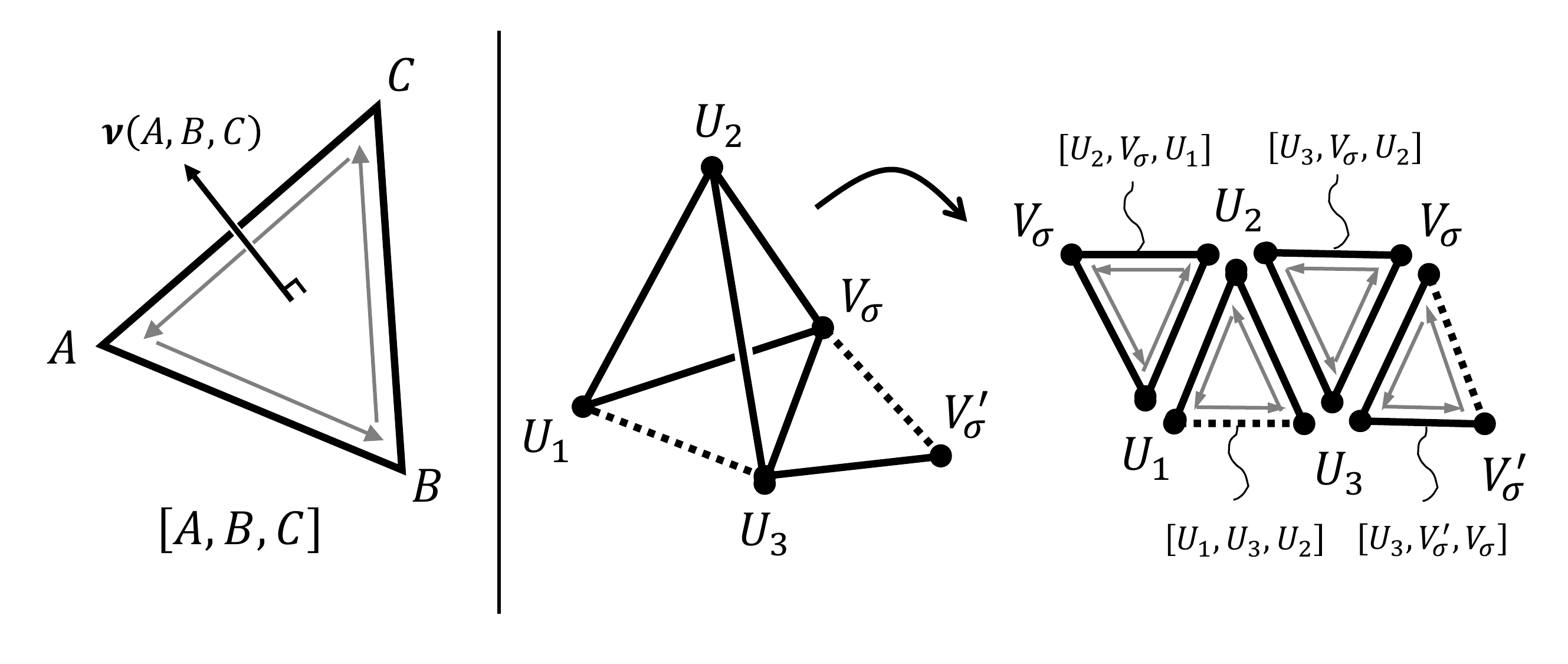}
              \caption{
                Oriented triangular faces and the fundamental patch $F_0$.
                Left: the boundary direction of $[A,B,C]$ and its unit
                normal are related by the right-hand rule.
                Center: the four triangles forming $F_0$.
                Right: the triangles shown separately, ordered from
                left to right as
                $[U_2,V_\sigma,U_1]$, $[U_1,U_3,U_2]$,
                $[U_3,V_\sigma,U_2]$, and $[U_3,V_\sigma',V_\sigma]$.
                Gray arrows indicate boundary directions, which are
                opposite along shared edges.
                Dotted edges $U_1U_3$ and $V_\sigma V_\sigma'$
                attach to the input and output caps, respectively.
              }
              \label{fig:oriented-triangle}
            \end{figure}
            
            The fundamental patch consists of four oriented triangles:
            \begin{equation}
              \begin{aligned}
                F_0 := {}&
                [U_1,U_3,U_2]+[U_2,V_\sigma,U_1]
                +[U_3,V_\sigma,U_2]+[U_3,V_\sigma',V_\sigma].
              \end{aligned}
              \label{eq:fundamental-face-set}
            \end{equation}
            Sums of faces represent gluing along edges with matching
            vertex labels, and their boundaries are computed by
            summing the oriented face boundaries.
            The second and third triangles in
            Eq.~\eqref{eq:fundamental-face-set} are the origami faces;
            the first and fourth are auxiliary faces spanning
            the input and output zigzags, respectively
            (Fig.~\ref{fig:oriented-triangle}).
            Contributions from shared edges cancel, leaving
            \begin{equation*}
              \partial F_0
              =
              [U_1,U_3]+[U_3,V_\sigma']
              +[V_\sigma',V_\sigma]+[V_\sigma,U_1].
            \end{equation*}
            
            Rotating this patch around the $X$-axis gives
            \begin{equation*}
              F:=\sum_{k=0}^{N-1}R^kF_0,
            \end{equation*}
            where $R$ acts on each face by rotating its vertices
            while preserving their ordering.
            Since $RU_1=U_3$ and $RV_\sigma=V_\sigma'$,
            the edges between neighboring rotational sectors
            cancel in the boundary sum, giving
            \begin{equation*}
              \partial F
              =
              \sum_{k=0}^{N-1}
              \left(
                [R^kU_1,R^{k+1}U_1]
                -
                [R^kV_\sigma,R^{k+1}V_\sigma]
              \right).
            \end{equation*}
            Here, rotational indices are understood modulo $N$.
            The two remaining boundary components are regular
            $N$-gons in planes perpendicular to the $X$-axis.
            
            We close them with the oriented polygonal caps
            \begin{equation*}
              \begin{aligned}
                D_I := [U_1,R^{N-1}U_1,\ldots,RU_1], \ \ 
                D_{I'} := [V_\sigma,RV_\sigma,\ldots,R^{N-1}V_\sigma].
              \end{aligned}
            \end{equation*}
            Their orientations cancel the corresponding boundaries
            of $F$ and give unit normals $-\mathbf{e}_X$ and
            $\mathbf{e}_X$, respectively.
            The orientations of the patch and caps are fixed by
            their vertex orderings, regardless of whether the output
            zigzag lies in the positive or negative axial direction
            relative to the input.
            
\begin{figure}[tb]
              \centering
              \includegraphics[width=\textwidth]
                {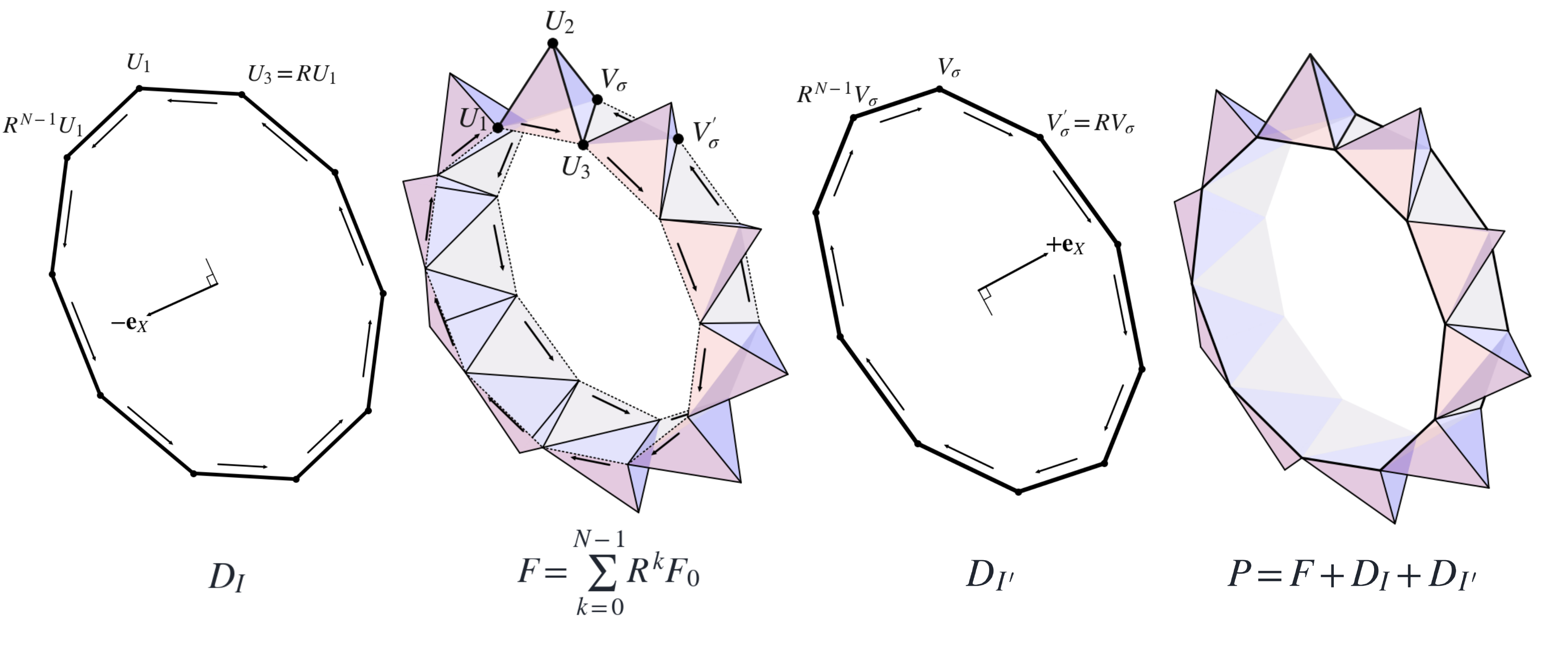}
              \caption{
                Construction of the closed oriented polyhedral surface $P$.
                From left to right: the input cap $D_I$;
                the surface $F$ formed from $N$ rotated copies of $F_0$;
                the output cap $D_{I'}$;
                and the assembled surface $P=F+D_I+D_{I'}$.
                Boundary arrows show how the cap boundaries cancel
                those of $F$.
                The input and output cap normals are
                $-\mathbf{e}_X$ and $\mathbf{e}_X$, respectively.
                Faces are shown with transparency to reveal the geometry.
              }
              \label{fig:closed-polyhedron}
            \end{figure}

            The assembled surface is
            \begin{equation*}
              P=F+D_I+D_{I'},
              \qquad
              \partial P=0.
            \end{equation*}
            Every edge belongs to exactly two faces with opposite
            induced boundary directions.
            Thus, the prescribed gluing defines a closed oriented
            polyhedral surface, possibly self-intersecting
            (Fig.~\ref{fig:closed-polyhedron}).
            
            With the caps treated as polygonal faces,
            the edges of $P$ form seven rotational classes,
            with representative lengths
            \begin{equation*}
              \begin{aligned}
                |U_1U_2| &= l_L, &
                |U_2U_3| &= l_R,\\
                |U_1V_\sigma| &= r_1, &
                |U_2V_\sigma| &= r_2, &
                |U_3V_\sigma| &= r_3,\\
                |U_1U_3| &= 2I, &
                |V_\sigma V_\sigma'| &= 2I'.
              \end{aligned}
            \end{equation*}
            Each class contains $N$ edges, counted according
            to the prescribed gluing.
            The first five lengths are fixed by the crease geometry,
            whereas the last two depend on the input and output states.

             \subsection{Derivation of the curvature action}
            \label{subsec:curvature-action}
            
            We apply the Euclidean Schl\"afli formula to the
            polyhedral surface $P$ constructed in
            Section~\ref{subsec:polyhedron}.
            Fix $N,l_L,l_R,r_1,r_2,r_3$ and the mountain--valley
            assignment specified by $\sigma$.
            All calculations are local, on a neighborhood where
            $f$ and the face normals are smooth.
            The angular variables and signed dihedral angles
            are represented by smooth real-valued lifts.
            
            Choose the lifts to agree on edges related by rotation.
            Write the five fixed-length-edge angles as
            $\Phi_{l_L}(\theta,I),\Phi_{l_R}(\theta,I),
            \Phi_{r_1}(\theta,I),\Phi_{r_2}(\theta,I),\Phi_{r_3}(\theta,I)$,
            and the angles at the edges of lengths $2I,2I'$ as
            $\Phi_I(\theta,I),\Phi_{I'}(\theta,I)$.
            The subscripts identify edge classes, not fixed angle values.
            All these functions refer to the selected folding solution
            with the prescribed geometric parameters fixed.
            We suppress their arguments in intermediate differential formulas.
            Since each class contains $N$ edges,
            the Schl\"afli formula, divided by $N$, reads
            \begin{equation}
              \begin{aligned}
                0={}&
                l_L\,\dd\Phi_{l_L}
                +l_R\,\dd\Phi_{l_R}
                +r_1\,\dd\Phi_{r_1}
                +r_2\,\dd\Phi_{r_2}
                +r_3\,\dd\Phi_{r_3}
                +2I\,\dd\Phi_I
                +2I'\,\dd\Phi_{I'}.
              \end{aligned}
              \label{eq:schlafli-per-sector}
            \end{equation}
            
            We first compute the last two terms.
            The input auxiliary face has unit normal
            \begin{equation*}
              \boldsymbol{\nu}(U_1,U_3,U_2)
              =
              -\sin\theta\,\mathbf{e}_X
              +\cos\theta\,\mathbf{e}_Z.
            \end{equation*}
            The boundary orientation of this face gives the unit
            tangent $-\mathbf{e}_Y$ along $[U_1,U_3]$,
            while the cap $D_I$ has normal $-\mathbf{e}_X$.
            Eq.~\eqref{eq:signed-dihedral-angle} therefore gives
            $\Phi_I\equiv\pi/2-\theta\pmod{2\pi}$.
            
            For the output segment, let
            $(\mathbf{e}_X,\mathbf{e}_{Y'},\mathbf{e}_{Z'})$
            be the right-handed frame defining $(\theta',I')$.
            The output auxiliary face has unit normal
            \begin{equation*}
              \boldsymbol{\nu}(U_3,V_\sigma',V_\sigma)
              =
              \sin\theta'\,\mathbf{e}_X
              -\cos\theta'\,\mathbf{e}_{Z'}.
            \end{equation*}
            Its boundary orientation gives the unit tangent
            $\mathbf{e}_{Y'}$ along $[V_\sigma',V_\sigma]$,
            while the cap $D_{I'}$ has normal $\mathbf{e}_X$.
            Thus, $\Phi_{I'}\equiv\theta'-\pi/2\pmod{2\pi}$.
            It follows that
            \begin{equation}
              \dd\Phi_I=-\dd\theta,
              \qquad
              \dd\Phi_{I'}=\dd\theta'.
              \label{eq:boundary-angle-differentials}
            \end{equation}
            
            Define
            \begin{equation}
              S_f(\theta,I)
              :=
              \frac12\left(
                l_L\Phi_{l_L}(\theta,I)
                +l_R\Phi_{l_R}(\theta,I)
                +r_1\Phi_{r_1}(\theta,I)
                +r_2\Phi_{r_2}(\theta,I)
                +r_3\Phi_{r_3}(\theta,I)
              \right),
              \label{eq:curvature-primitive}
            \end{equation}
            where the dihedral angles are functions of the input
            state $(\theta,I)$.
            The five edge lengths in this expression are fixed.
            Substituting Eq.~\eqref{eq:boundary-angle-differentials}
            into Eq.~\eqref{eq:schlafli-per-sector} therefore yields
            \begin{equation}
              \dd S_f
              =
              I\,\dd\theta-I'\,\dd\theta'.
              \label{eq:curvature-variation}
            \end{equation}
            Taking the exterior derivative gives
            \begin{equation}
              \dd\theta'\wedge\dd I'
              =
              \dd\theta\wedge\dd I.
              \label{eq:elementary-area-preservation}
            \end{equation}
            Thus, $f$ preserves the symplectic form
            $\dd\theta\wedge\dd I$.
            Neither the definition of $S_f$ nor
            Eq.~\eqref{eq:curvature-variation} requires a twist condition.
            
            Under the twist condition
            $\left.\partial\theta'/\partial I\right|_\theta\neq0$,
            we may write $I=I(\theta,\theta')$ locally.
            By Eq.~\eqref{eq:curvature-variation} and
            Section~\ref{subsec:variational-principle},
            $L_f(\theta,\theta')
            :=S_f\bigl(\theta,I(\theta,\theta')\bigr)$
            is a local generating function for $f$.
            Thus $S_f$ is the contribution of the fixed-length edges to the
            signed total discrete mean curvature of $P$~\cite{Sullivan2008Curvatures},
            divided by $N$, the number of rotational sectors.
            The edges of lengths $2I$ and $2I'$ are omitted from
            this action; their terms in the Schl\"afli formula
            give the right-hand side of
            Eq.~\eqref{eq:curvature-variation}.
            
            Any two smooth lifts of a given dihedral angle differ
            locally by a constant integer multiple of $2\pi$.
            Since the coefficients in Eq.~\eqref{eq:curvature-primitive}
            are fixed, a change of lifts alters $S_f$ and $L_f$
            only by an additive constant.
            The differential identity and the generating relations
            are consequently unchanged.

        \subsection{Reparameterization and composition of the curvature action}
        \label{subsec:composition}
        
\begin{figure}[tb]
          \centering
          \includegraphics[width=\textwidth]
            {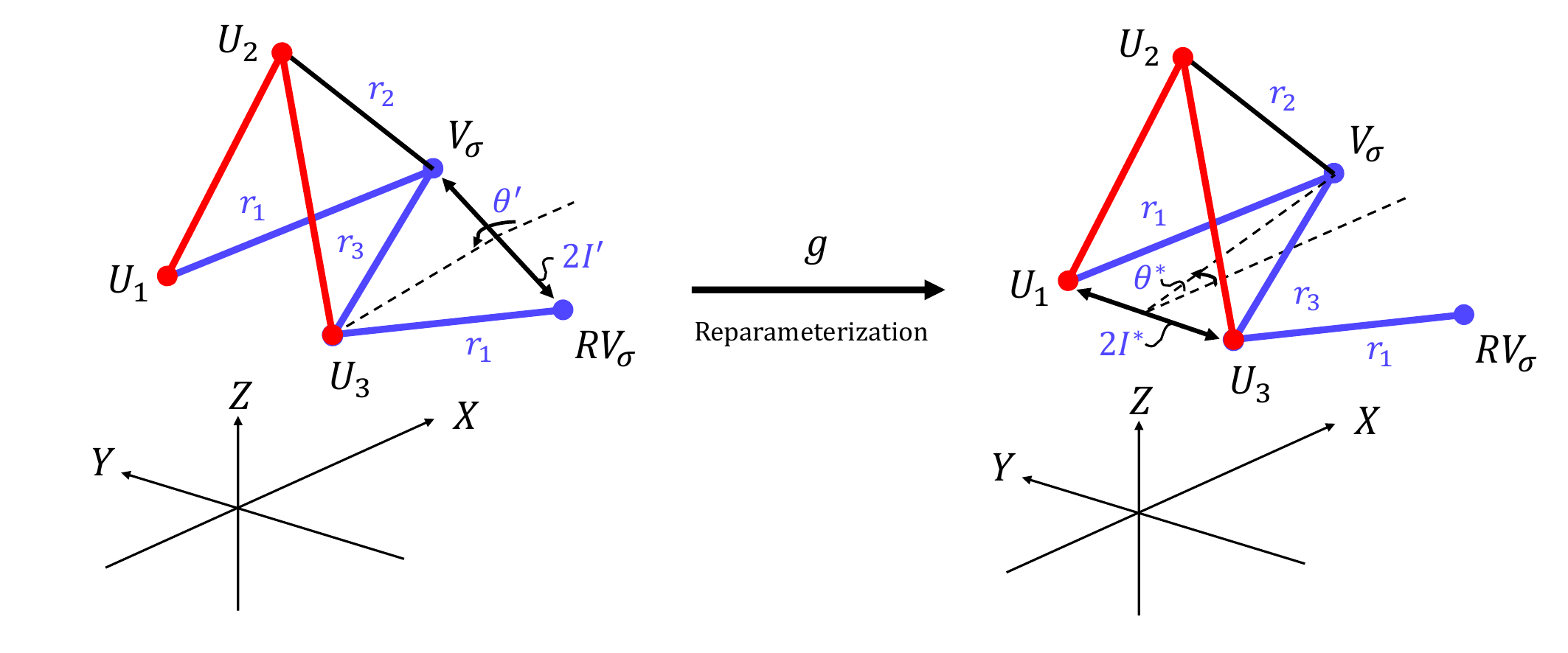}
          \caption{
            Reparameterization of the output zigzag.
            Left: $(V_\sigma,U_3,RV_\sigma)$ defines
            $(\theta',I')$, with $2I'=|V_\sigma RV_\sigma|$.
            Right: $(U_1,V_\sigma,U_3)$ defines
            $(\theta^*,I^*)$, with $2I^*=|U_1U_3|=2I$.
            Red and blue edges belong to the input and output
            zigzags of $f$, respectively.
            The vertex positions and crease lengths are unchanged.
          }
          \label{fig:reparameterization}
        \end{figure}

        The preceding construction gives a curvature action $S_f$
        for the elementary propagation map $f$.
        We now consider the case in which $g$ is applied after $f$.
        Write $\mathbf{x}'=f(\mathbf{x})$ and
        $\mathbf{x}^*=g(\mathbf{x}')$.
        As shown in Fig.~\ref{fig:reparameterization},
        $g$ replaces the representative segment
        $(V_\sigma,U_3,V_\sigma')$ by $(U_1,V_\sigma,U_3)$,
        where $V_\sigma'=RV_\sigma$.
        The folded zigzag is unchanged.
        Since the input segment and the final output segment
        have the same endpoints,
        $I^*=\frac12|U_1U_3|=I$.
        
        For $g\circ f$, the two origami faces and the input
        and output auxiliary faces form the oriented surface
        \begin{equation*}
          \begin{aligned}
            F_{g\circ f}:={}&
            [U_1,U_3,U_2]+[U_2,V_\sigma,U_1]
            +[U_3,V_\sigma,U_2]+[U_1,V_\sigma,U_3].
          \end{aligned}
        \end{equation*}
        Each edge occurs twice with opposite boundary directions,
        so $\partial F_{g\circ f}=0$.
        We may therefore apply the Euclidean Schl\"afli formula
        directly to $F_{g\circ f}$.
        Its five fixed edge lengths are $l_L,l_R,r_1,r_2,r_3$,
        and its remaining edge $U_1U_3$ has length $2I$.
        All faces are assumed nondegenerate, and the signed
        dihedral angles are given smooth local real-valued lifts.
        
        Denote these angles by
        $\Psi_{l_L}(I),\Psi_{l_R}(I),\Psi_{r_1}(I),
        \Psi_{r_2}(I),\Psi_{r_3}(I)$ and $\Psi_I(I)$, respectively.
        On the selected solution, the six edge lengths determine the
        tetrahedral shape. Thus these angles depend only on $I$;
        varying $\theta$ at fixed $I$ rotates the tetrahedron rigidly
        about $U_1U_3$. We suppress the arguments in intermediate formulas.
        The input and output auxiliary normals are
        $-\sin\theta\,\mathbf{e}_X+\cos\theta\,\mathbf{e}_Z$
        and
        $\sin\theta^*\,\mathbf{e}_X-\cos\theta^*\,\mathbf{e}_Z$.
        Using the boundary tangent $-\mathbf{e}_Y$
        induced by the input face gives
        $\Psi_I\equiv\theta^*-\theta+\pi\pmod{2\pi}$,
        and hence $\dd\Psi_I=\dd\theta^*-\dd\theta$.
        
        Define
        \begin{equation*}
          S_{g\circ f}(I)
          :=
          \frac12\left(
            l_L\Psi_{l_L}(I)
            +l_R\Psi_{l_R}(I)
            +r_1\Psi_{r_1}(I)
            +r_2\Psi_{r_2}(I)
            +r_3\Psi_{r_3}(I)
          \right).
        \end{equation*}
        Since these five lengths are fixed,
        the Schl\"afli formula gives
        \begin{equation}
          \dd S_{g\circ f}
          =
          -I\,\dd\Psi_I
          =
          I\,\dd\theta-I^*\,\dd\theta^*.
          \label{eq:gf-action-differential}
        \end{equation}
        This has the same form as
        Eq.~\eqref{eq:curvature-variation} for $S_f$.
        In both cases, the action is one half of the
        edge-length-weighted sum of signed dihedral angles
        at the five fixed-length edges per rotational sector.
        The auxiliary faces used to define those angles
        depend on whether the step is $f$ or $g\circ f$.
        
        For the steps $M_i$ defined in
        Section~\ref{subsec:zigzag-propagation}, set
        \begin{equation*}
          S_i(\theta,I)=
          \begin{cases}
            S_{f_i}(\theta,I), & M_i=f_i,\\
            S_{g_i\circ f_i}(I), & M_i=g_i\circ f_i.
          \end{cases}
        \end{equation*}
        For one iteration $\mathbf x_{n+1}=M(\mathbf x_n)$, introduce
        intermediate states
        $\mathbf x_{n,i}=(\theta_{n,i},I_{n,i})^{\mathsf T}$ satisfying
        \[
          \mathbf x_{n,0}=\mathbf x_n,\qquad
          \mathbf x_{n,i+1}=M_i(\mathbf x_{n,i}),\qquad
          0\leq i<m.
        \]
        Then $\mathbf x_{n,m}=\mathbf x_{n+1}$.
        Regarding the intermediate states as functions of $\mathbf x_n$,
        define
        \begin{equation}
          S_M(\mathbf x_n)
          :=\sum_{i=0}^{m-1}S_i(\theta_{n,i},I_{n,i}).
          \label{eq:composite-curvature-action}
        \end{equation}
        The intermediate terms cancel upon differentiation:
        \begin{equation}
          \begin{aligned}
            \dd S_M(\mathbf x_n)
            &=\sum_{i=0}^{m-1}
              \left(I_{n,i}\,\dd\theta_{n,i}
                    -I_{n,i+1}\,\dd\theta_{n,i+1}\right)\\
            &=I_n\,\dd\theta_n-I_{n+1}\,\dd\theta_{n+1}.
          \end{aligned}
          \label{eq:composite-action-differential}
        \end{equation}
        Under the twist condition for $M$, expressing $S_M(\mathbf x_n)$
        in terms of $(\theta_n,\theta_{n+1})$ gives a local generating
        function $L(\theta_n,\theta_{n+1})$ for $M$.
        The discrete variational principle of
        Section~\ref{subsec:variational-principle} therefore
        applies to this curvature action.

        To interpret the sum geometrically, consider the
        auxiliary triangles at an intermediate zigzag.
        After identifying their representatives by rotation,
        the output auxiliary triangles of one step coincide
        with the input auxiliary triangles of the next,
        with opposite orientations.
        Removing them joins the adjacent origami faces.
        
        The effect on the dihedral angles can be computed
        along each edge of the common zigzag.
        Let $\mathbf{a}$ and $\mathbf{b}$ be the oriented
        unit normals of the two origami faces, and let
        $\mathbf{n}$ and $-\mathbf{n}$ be those of the
        auxiliary faces to be removed.
        Choose the edge tangent $\mathbf{t}$ according
        to the boundary orientation of the first origami face.
        Writing $\angle_{\mathbf{t}}$ for the signed angle
        about $\mathbf{t}$, the two contributions are
        $\Phi_e^-=\angle_{\mathbf{t}}(\mathbf{a},\mathbf{n})$
        and
        $\Phi_e^+=\angle_{\mathbf{t}}(-\mathbf{n},\mathbf{b})$.
        The angle between the remaining faces is
        $\varphi_e=\angle_{\mathbf{t}}(\mathbf{a},\mathbf{b})$.
        Hence,
        \begin{equation*}
          \Phi_e^-+\Phi_e^+
          =
          \varphi_e+\pi+2\pi k_e,
          \qquad k_e\in\mathbb{Z},
        \end{equation*}
        where $k_e$ is locally constant for smooth angle lifts.
        Since $\ell_e$ is fixed, the last two terms contribute
        only an additive constant to the action.
        
        It follows that
        \begin{equation}
          S_M(\mathbf x_n)
          =
          \frac12
          \left(
            \sum_{e\in E_{\mathrm{int}}}\ell_e\varphi_e(\mathbf x_n)
            +
            \sum_{e\in E_{\partial}}\ell_e\beta_e(\mathbf x_n)
          \right)
          +C,
          \label{eq:composite-crease-curvature}
        \end{equation}
        where each sum contains one representative edge
        from each rotational class.
        The set $E_{\mathrm{int}}$ consists of the creases
        shared by two origami faces, with signed dihedral
        angles $\varphi_e$.
        The set $E_{\partial}$ consists of the edges of
        the two end zigzags, where $\beta_e$ is the angle
        between the adjacent origami face and the retained
        auxiliary face.
        The angles $\varphi_e(\mathbf x_n)$ and $\beta_e(\mathbf x_n)$
        depend on the input state through the selected propagation.
        All edge lengths in both sums are fixed, and $C$
        is independent of the state on the chosen neighborhood.
        
        Thus, the dependence on internal auxiliary faces
        is eliminated, while the two ends retain contributions
        from the auxiliary faces used to close the strip.
        These boundary contributions are part of the curvature
        action and are retained in the variational formulation;
        fixing the endpoint angular variables does not,
        in general, make them constant.

    \subsection{Discrete curvature action for conformally symplectic propagation}
        \label{sec:conformal}
        
        A map $F$ is conformally symplectic if $F^*\omega=\lambda\omega$
        for a constant $\lambda>0$. The dissipative standard map is a
        representative example with $0<\lambda<1$~\cite{Calleja2022Dissipative}.
        Under a twist condition, orbits of a conformally symplectic map can
        be characterized locally by the stationarity of a discrete action
        in which the generating function at step $n$ is weighted by
        $\lambda^{-n}$.
        The tubular origami considered here also gives rise to a
        conformally symplectic system when the scale changes from one
        ring to the next~\cite{Imada2023Undulations,Ichikawa2026InvariantCurves}.
        We now derive the discrete curvature action in this case using
        the same geometric construction as in the preceding subsection
        for the conservative case.

        Let the crease lengths of the output zigzag be $sl_L$ and $sl_R$,
        where $s>0$. We normalize the output by dividing its length
        coordinate by $s$, and write the normalized propagation as
        \[
          \widehat M(\theta_n,I_n)=(\theta_{n+1},I_{n+1}).
        \]
        Thus the unnormalized output half-chord length is $sI_{n+1}$.
        By Eq.~\eqref{eq:input-domain}, both normalized states satisfy
        \[
          (\theta_k,I_k)\in
          \mathbb R/(2\pi\mathbb Z)\times
          \left(\frac{|l_L-l_R|}{2},\,\frac{l_L+l_R}{2}\right),
          \qquad k=n,n+1.
        \]
        The boundary variation in
        Eq.~\eqref{eq:composite-action-differential} also applies when
        the input and output crease lengths differ, since its derivation
        uses fixed crease lengths within each elementary step.
        Substituting the unnormalized output length gives
        \[
          \dd S_M(\theta_n,I_n)
          =I_n\,\dd\theta_n-sI_{n+1}\,\dd\theta_{n+1}.
        \]
        Taking the exterior derivative yields
        \begin{equation}
          \widehat M^*\omega=s^{-1}\omega,\qquad
          \omega=\dd\theta\wedge\dd I,
          \label{eq:tube-conformal-factor}
        \end{equation}
        and hence $\det D\widehat M=s^{-1}$.
        Under the twist condition, express the reference-module action as
        $L(\theta_n,\theta_{n+1})$ to obtain
        \[
          I_n=D_1L(\theta_n,\theta_{n+1}),\qquad
          I_{n+1}=-s^{-1}D_2L(\theta_n,\theta_{n+1}).
        \]

        If the crease pattern at step $n$ is scaled by $s^n$ relative to the
        reference pattern, each curvature term, being a product of an edge
        length and a dihedral angle, scales by $s^n$. The discrete curvature
        action is therefore
        \begin{equation}
          \mathcal A(\theta_0,\ldots,\theta_J)
          =\sum_{n=0}^{J-1}s^nL(\theta_n,\theta_{n+1}).
          \label{eq:tube-weighted-action}
        \end{equation}
        With the endpoint angles fixed, stationarity at each interior ring
        gives
        \begin{equation}
          D_2L(\theta_{n-1},\theta_n)
          +sD_1L(\theta_n,\theta_{n+1})=0,
          \qquad 1\leq n<J.
          \label{eq:tube-weighted-euler-lagrange}
        \end{equation}
        Thus the orbits of $\widehat M$ are locally characterized by
        the stationary points of $\mathcal A$.

\section{Symplectic propagation in the Maxwell origami tube}
    \label{sec:maxwell}
    
    In the Maxwell origami tube of Imada and Tachi~\cite{Imada2025Maxwell},
    ring-to-ring propagation defines a map on a boundary state space of
    dimension \(4N-6\), where \(N\) is the number of circumferential units. 
    When the crease data and mountain--valley assignments are invariant
    under circumferential shifts, the \(N\)-fold rotationally symmetric
    configurations form a two-dimensional invariant subspace. The restriction of the propagation
    map to this subspace is a particular instance of the tubular origami maps
    considered in the preceding section. 
    
    We introduce an alternative to the boundary coordinates used by Imada and Tachi. For each diagonal of a triangulated boundary polygon, we take its half-length and the signed dihedral angle between the two adjacent triangles as a pair of coordinates.
    The discrete curvature action then shows that the smooth local propagation map determined by the prescribed mountain--valley assignment is symplectic, even when the left and right edges of the zigzag have unequal prescribed lengths.
    
    \subsection{Propagation between boundary polygons}
        \label{subsec:maxwell_propagation}
            \begin{figure}[tb]
              \centering
              \includegraphics[width=\textwidth]
                {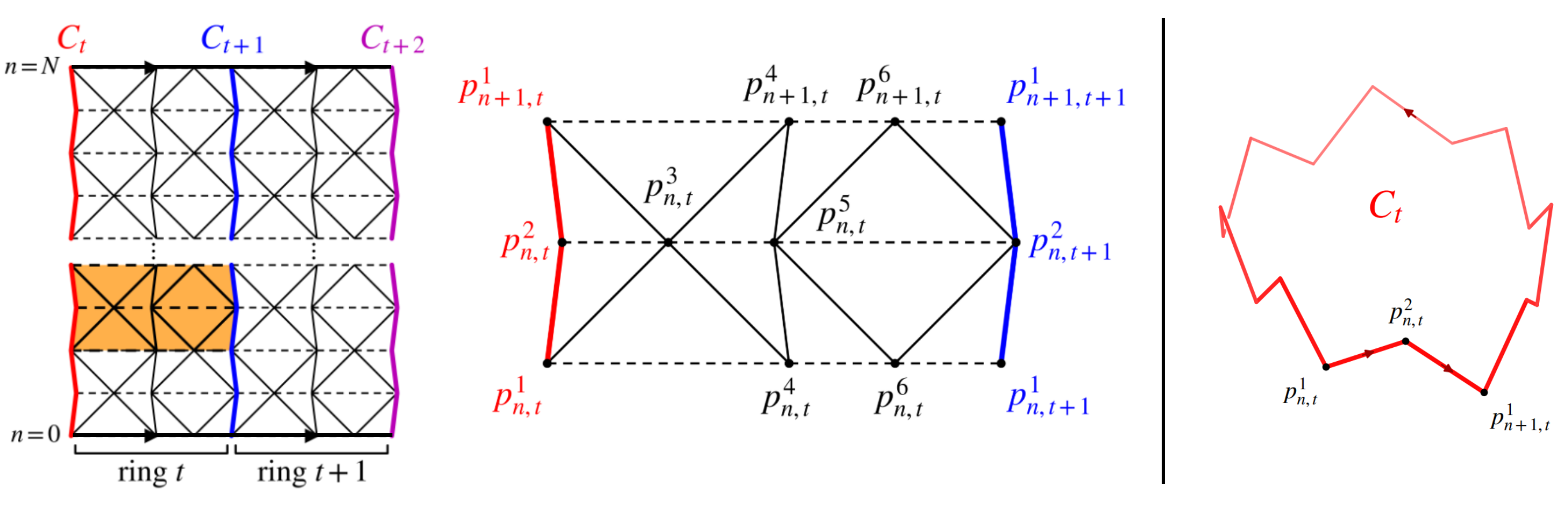}
                \caption{Crease pattern and oriented zigzag boundary of the Maxwell origami tube.
                In the crease pattern (left), the edges at $n=0$ and $n=N$ are identified
                as indicated by the black arrows. The red, blue, and purple zigzags
                correspond to the successive boundaries $C_t$, $C_{t+1}$, and $C_{t+2}$,
                with $C_{t+1}$ shared by rings $t$ and $t+1$.
                The center panel shows the cell $(n,t)$ highlighted in orange on the left.
                Here $p^i_{n,t}\in\mathbb{R}^3$ denotes the position of vertex $i$ in
                the folded configuration. Each cell contributes the oriented edges
                $[p^1_{n,t},p^2_{n,t}]$ and $[p^2_{n,t},p^1_{n+1,t}]$ to $C_t$,
                with $n$ taken modulo $N$. Since no rotational symmetry is imposed,
                the cells around a ring may have different folded configurations.
                The right panel shows a nonplanar realization of $C_t$ with $N=8$
                and no rotational symmetry. The red arrows indicate its orientation.}
              \label{fig:maxwell_model}
            \end{figure}
            
            We use the crease connectivity of the Maxwell origami tube shown in
            Fig.~\ref{fig:maxwell_model}. Let $N\geq 3$ be the number of cells
            around the tube, and let $p^i_{n,t}\in\mathbb{R}^3$ denote the
            position of vertex $i$ of cell $(n,t)$ after folding. Taking the
            circumferential index $n$ modulo $N$, we write the input boundary
            of ring $t$ as the closed oriented polygonal chain
            \begin{equation*}
              C_t :=
              \sum_{n=0}^{N-1}
              \bigl(
                [p^1_{n,t},p^2_{n,t}]
                +[p^2_{n,t},p^1_{n+1,t}]
              \bigr).
            \end{equation*}
            Following the orientation of $C_t$, denote the prescribed lengths
            of the two edges in cell $(n,t)$ by
            $l_{\mathrm L,n},l_{\mathrm R,n}\in\mathbb{R}_{>0}$, respectively:
            \begin{equation*}
              \lVert p^2_{n,t}-p^1_{n,t}\rVert=l_{\mathrm L,n},
              \qquad
              \lVert p^1_{n+1,t}-p^2_{n,t}\rVert=l_{\mathrm R,n}.
            \end{equation*}
            These lengths are independent of $t$ because the crease pattern
            repeats from one ring to the next. In the reference pattern of
            Imada and Tachi~\cite{Imada2025Maxwell}, mirror symmetry gives
            $l_{\mathrm L,n}=l_{\mathrm R,n}$. Here we allow the two lengths
            to differ.
            
            For each $n$, the chord joining $p^1_{n,t}$ to $p^1_{n+1,t}$ has
            length
            \[
              d^1_{n,t}:=\lVert p^1_{n+1,t}-p^1_{n,t}\rVert .
            \]
            This chord and the two boundary edges form a nondegenerate triangle
            precisely when
            \[
              |l_{\mathrm L,n}-l_{\mathrm R,n}|
              <d^1_{n,t}<
              l_{\mathrm L,n}+l_{\mathrm R,n}.
            \]
            The chord length depends on the folded configuration. Even when
            $l_{\mathrm L,n}=l_L$ and $l_{\mathrm R,n}=l_R$ for all $n$,
            $d^1_{n,t}$ may vary around the ring, since no rotational symmetry
            is imposed on the folded tube.
            
            To construct the next boundary, consider first the vertex $p^3_{n,t}$
            in the central panel of Fig.~\ref{fig:maxwell_model}. It is connected
            by creases to $p^1_{n,t}$, $p^2_{n,t}$, and $p^1_{n+1,t}$, whose
            positions are given by $C_t$. Hence $p^3_{n,t}$ lies at the
            intersection of three spheres centered at these vertices, with radii
            equal to the lengths of the corresponding creases. The prescribed
            mountain--valley assignment selects the intersection point. This
            construction also applies when $l_{\mathrm L,n}\ne l_{\mathrm R,n}$.
            
            Applying the same construction around the ring gives the successive
            sets of vertices
            \[
              p^3_{n,t},\quad p^4_{n,t},\quad p^5_{n,t},\quad p^6_{n,t},
              \quad p^2_{n,t+1},\quad p^1_{n,t+1}.
            \]
            The last two sets form the output boundary $C_{t+1}$.
            Where all six stages admit transverse sphere intersections and the
            mountain--valley assignment determines a unique solution at each
            stage, the construction defines a smooth local propagation map
            $M$ from $C_t$ to $C_{t+1}$. The explicit vertex construction for
            the original crease pattern is given in the Supplemental Material
            of Imada and Tachi~\cite{Imada2025Maxwell}.

    \subsection{Diagonal lengths and dihedral angles as boundary coordinates}
        \label{subsec:maxwell_boundary_coordinates}

        Boundary shapes that differ only by translation and rotation are
        identified, with the cyclic vertex labels fixed.
        
        Imada and Tachi reconstruct the vertices of a boundary zigzag
        successively around the circumference. In their coordinates,
        $d^i_{n,t}$ measures a chord spanning two consecutive zigzag edges,
        whereas $\rho^i_{n,t}$ measures the dihedral angle between successive
        auxiliary triangles about their shared zigzag edge. A further angle
        $\eta_t$ specifies the position of the last vertex.
        In the two-dimensional system of Section~\ref{sec:action}, the half-length $I$ and
        the angular differential $\dd\theta$ arise from the same auxiliary edge
        in the Schl\"afli formula. We extend this pairing to the Maxwell
        boundary by associating each diagonal of a triangulation with the
        dihedral angle measured about that diagonal.
        
        The prescribed lengths of the boundary edges of $C_t$ still allow
        its vertices to move relative to one another. To see what additional
        data determine its shape, consider a spatial quadrilateral with four
        prescribed side lengths. Choose one diagonal. Once its length is
        given, all three sides of each of the two resulting triangles are
        known. The two triangles can still rotate about their common
        diagonal. Their signed dihedral angle specifies this relative
        rotation and hence determines the shape of the quadrilateral.
        The same construction applies to a polygon with more vertices.
        Take a convex reference polygon with the same cyclic ordering of
        vertices and draw diagonals until it is divided into triangles,
        with no two diagonals crossing in its interior. This choice of
        triangles is a triangulation. Joining the corresponding vertices
        in space gives a triangulated surface spanning the boundary, even
        when its vertices are noncoplanar. The diagonal lengths determine
        the shapes of the triangles, and the signed dihedral angles
        determine how adjacent triangles fit together.
        
        This geometric description also gives canonical coordinates.
        Kapovich and Millson showed that rotation about a diagonal is the
        Hamiltonian motion generated by its length. For a fan triangulation
        with nondegenerate triangles, they proved that the diagonal lengths
        and the signed bending angles at the same diagonals form local
        action--angle coordinates
        \cite{kapovich1996polygons}.
        More generally, the bending systems of
        Klyachko~\cite{klyachko1994spatial} and
        Kapovich and Millson~\cite{kapovich1996polygons} associate a
        commuting family of action variables to the diagonal lengths of
        any triangulation; see also
        \cite{nohara2014toric}.
        
        To construct a cap for $C_t$, first join the two neighbors of each
        vertex $p^2_{n,t}$. This gives one triangle for each zigzag segment
        and leaves a central $N$-gon with vertices
        $p^1_{0,t},\ldots,p^1_{N-1,t}$. In the rotationally symmetric
        setting of Section~\ref{sec:action}, these vertices lie in a common plane, so the
        central polygon itself serves as a single face of the cap. In the
        Maxwell tube, the vertices may be noncoplanar. We therefore span
        the central polygon by a fan of triangles, joining $p^1_{0,t}$ to
        $p^1_{k,t}$ for $2\leq k\leq N-2$. Together with the $N$ triangles
        along the zigzag, this fan gives a triangulated cap with oriented
        boundary $C_t$.
        We use diagonal half-lengths rather than full lengths to match the
        factor $1/2$ in the curvature action. The signs of the angular
        variables follow the orientation convention in
        Eq.~\eqref{eq:signed-dihedral-angle}.
        The cap has $N+(N-3)=2N-3$ interior edges. For each such edge
        $e_{j,t}$, we take its half-length $I_{j,t}=|e_{j,t}|/2$ and its
        signed exterior dihedral angle $\theta_{j,t}$ as a pair of boundary
        variables.
        Where all triangles are nondegenerate, these $2N-3$ pairs thus give
        $4N-6$ local coordinates for the boundary shape, in agreement with
        the number of degrees of freedom of one Maxwell ring. We use these coordinates below to express
        the curvature action for propagation from $C_t$ to $C_{t+1}$.

    \subsection{Oriented caps and the curvature action of one ring}
    \label{subsec:maxwell-ring-action}
    
            \begin{figure}[tb]
              \centering
              \includegraphics[width=\textwidth]
                {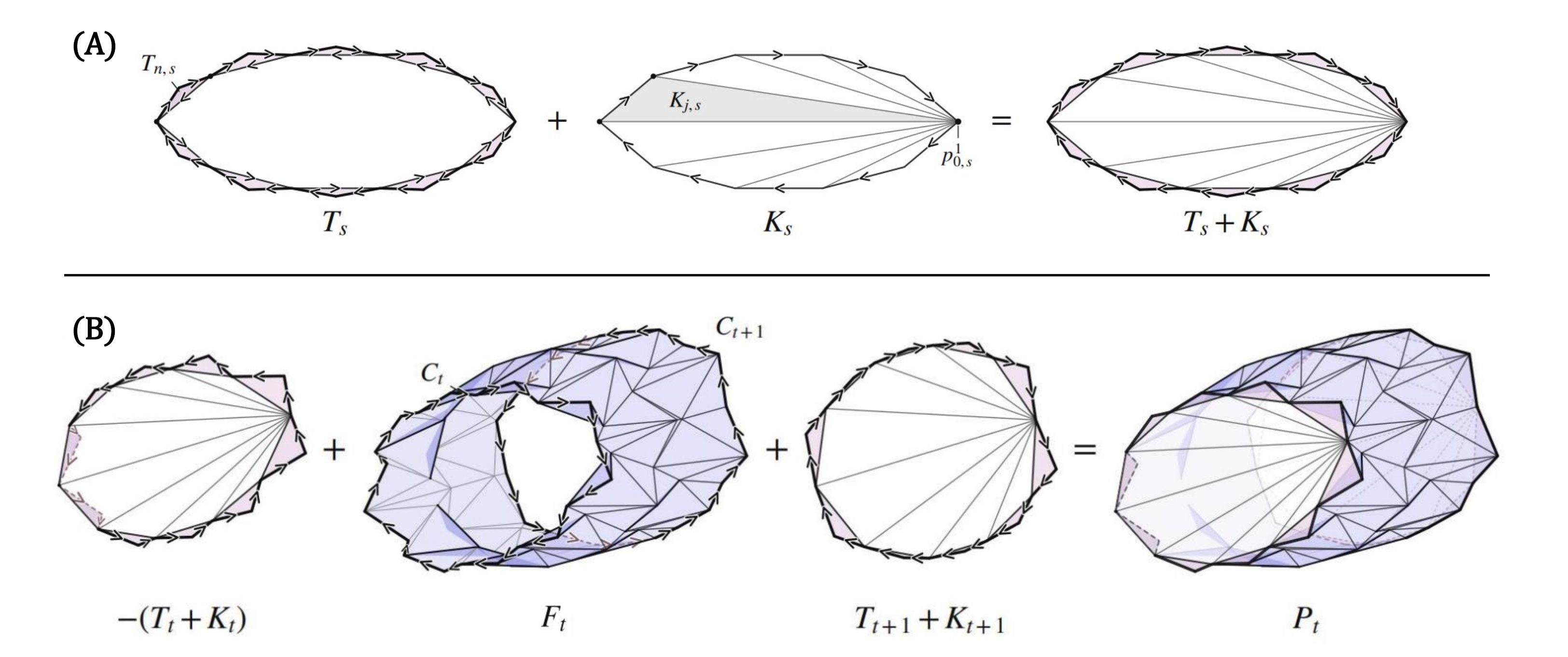}
                \caption{Oriented caps and closure of a single Maxwell origami ring.
                (A) At each boundary zigzag $C_s$, the triangles $T_{n,s}$ form $T_s$,
                and the triangles $K_{j,s}$ form a fan $K_s$ from $p^1_{0,s}$
                spanning the polygon with vertices $p^1_{n,s}$.
                The arrows indicate the induced boundary orientations.
                Along each chord $[p^1_{n,s},p^1_{n+1,s}]$, the two orientations
                cancel, leaving $\partial(T_s+K_s)=C_s$.
                (B) The origami surface $F_t$ has oriented boundary $C_t-C_{t+1}$.
                Attaching the input cap $-(T_t+K_t)$ along $C_t$ and the output cap
                $T_{t+1}+K_{t+1}$ along $C_{t+1}$ yields the closed oriented
                polyhedral surface $P_t$.}
              \label{fig:maxwell-caps}
            \end{figure}
    
        We use the triangulation of the preceding subsection to close the
        origami surface between $C_t$ and $C_{t+1}$. For either boundary,
        $s\in\{t,t+1\}$, orient the auxiliary triangles by
        \[
          T_{n,s}=[p^1_{n,s},p^2_{n,s},p^1_{n+1,s}],
          \qquad
          K_{j,s}=[p^1_{0,s},p^1_{j,s},p^1_{j+1,s}],
        \]
        where $n$ is taken modulo $N$ and $1\leq j\leq N-2$. Write
        \begin{equation*}
          T_s=\sum_{n=0}^{N-1}T_{n,s},
          \qquad
          K_s=\sum_{j=1}^{N-2}K_{j,s}.
        \end{equation*}
        The triangles of $T_s$ follow the zigzag, while those of $K_s$ span
        the central polygon. Along the chord $[p^1_{n,s},p^1_{n+1,s}]$,
        the central triangulation $K_s$ induces the displayed direction
        and $T_{n,s}$ induces the opposite direction. Thus
        $\partial(T_s+K_s)=C_s$. This construction also applies when the
        vertices of the central polygon are noncoplanar, since each triangle
        carries its own oriented normal.

        Orient the origami surface $F_t$ of ring $t$ so that
        $\partial F_t=C_t-C_{t+1}$. Attaching the input and output caps gives
        the closed oriented polyhedral surface
        \[
          P_t=F_t-(T_t+K_t)+(T_{t+1}+K_{t+1}),\qquad \partial P_t=0.
        \]
        
        For either boundary $s\in\{t,t+1\}$, index the $q=2N-3$
        interior diagonals of $T_s+K_s$ by $e_{j,s}$, $1\leq j\leq q$.
        Set $I_{j,s}=|e_{j,s}|/2$ and let $\theta_{j,s}$ be the signed
        exterior dihedral angle at $e_{j,s}$, using
        Eq.~\eqref{eq:signed-dihedral-angle}. Along $[p^1_{n,s},p^1_{n+1,s}]$,
        this angle is measured between $T_{n,s}$ and the incident triangle
        of $K_s$. Along $[p^1_{0,s},p^1_{k,s}]$, $2\leq k\leq N-2$, it is
        measured between $K_{k-1,s}$ and $K_{k,s}$.

        Using locally smooth real-valued lifts of the angles, form the column
        vectors
        \begin{equation*}
          \boldsymbol{\theta}_s
            =(\theta_{1,s},\ldots,\theta_{q,s})^{\mathsf T},
          \qquad
          \mathbf I_s
            =(I_{1,s},\ldots,I_{q,s})^{\mathsf T}.
        \end{equation*}
        In these local coordinates, the boundary state and its propagation
        are given by
        \begin{equation}
          \mathbf x_t=
          \begin{pmatrix}
            \boldsymbol{\theta}_t\\
            \mathbf I_t
          \end{pmatrix}
          \in\mathbb R^{2q},
          \qquad
          \mathbf x_{t+1}=M(\mathbf x_t).
          \label{eq:maxwell-state-map}
        \end{equation}

        Let $E_{\mathrm{fix}}$ be the set of origami edges, including
        the sides of both zigzags. Write $\Phi_e(\mathbf x_t)$ for the
        signed exterior dihedral angle at $e$ in $P_t$, with the prescribed
        crease lengths and mountain--valley assignment fixed. At a zigzag
        side, the angle is measured between an origami face and a cap
        triangle. The curvature action of one ring is
        \[
          S_M(\mathbf x_t)
          =\frac12\sum_{e\in E_{\mathrm{fix}}}\ell_e\Phi_e(\mathbf x_t).
        \]
        In $P_t$, the orientation of the input cap is reversed and that of
        the output cap is preserved. Thus their angle differentials are
        $-\dd\theta_{j,t}$ and $\dd\theta_{j,t+1}$, respectively.
        Since the origami edge lengths are fixed, the Euclidean
        Schl\"afli formula on $P_t$ gives
        \begin{equation}
          \dd S_M(\mathbf x_t)
          =\mathbf I_t^{\mathsf T}\dd\boldsymbol{\theta}_t
           -\mathbf I_{t+1}^{\mathsf T}\dd\boldsymbol{\theta}_{t+1}.
          \label{eq:maxwell-vector-action}
        \end{equation}
        This identity holds for unequal prescribed left and right edge lengths.
        On a common regular domain, changing the cap triangulations changes the
        boundary coordinates and adds only separate input and output boundary
        terms to the one-ring curvature action, so the resulting variational
        principles describe the same propagation between boundary shapes.
        Taking the exterior derivative of Eq.~\eqref{eq:maxwell-vector-action}
        gives
        \begin{equation}
          M^*\omega=\omega,
          \qquad
          \omega=\sum_{j=1}^{q}\dd\theta_j\wedge\dd I_j.
          \label{eq:maxwell-symplectic-form}
        \end{equation}
        Eq.~\eqref{eq:maxwell-symplectic-form} also follows from Miranda's
        isotropy theorem applied to the two oriented boundaries of the
        ring~\cite{Miranda2025Domes}. The input and output boundaries carry
        opposite orientations, so isotropy on the graph of $M$ gives
        $\omega_{\mathrm{in}}-M^*\omega_{\mathrm{out}}=0$. This preservation
        law implies that the boundary one-form in
        Eq.~\eqref{eq:maxwell-vector-action} has a local primitive. At a fixed
        point on a smooth branch of the propagation map, $DM$ is symplectic,
        and its eigenvalues occur in reciprocal pairs counted with
        algebraic multiplicity.

        The Schl\"afli formula, however, allows us to construct this primitive
        directly from the folded geometry without integrating the propagation
        relations. Eq.~\eqref{eq:maxwell-vector-action} identifies the
        primitive as $S_M$, the curvature sum over the fixed-length creases
        and boundary edges of the capped ring.
        
        Under the twist condition
        \begin{equation*}
          \det\!\left(
            \frac{\partial\boldsymbol{\theta}_{t+1}}{\partial\mathbf I_t}
          \right)\ne 0,
        \end{equation*}
        expressing this same curvature sum in terms of the input and output
        boundary angles gives the local generating function
        $L(\boldsymbol{\theta}_t,\boldsymbol{\theta}_{t+1})$, with generating relations
        \[
          \mathbf I_t=D_1L(\boldsymbol{\theta}_t,\boldsymbol{\theta}_{t+1}),
          \qquad
          \mathbf I_{t+1}=-D_2L(\boldsymbol{\theta}_t,\boldsymbol{\theta}_{t+1}).
        \]
        With the endpoint vectors $\boldsymbol{\theta}_0$ and
        $\boldsymbol{\theta}_J$ fixed, the angle sequences of propagation
        orbits are precisely the stationary points of the discrete action
        \begin{equation}
          \mathcal A(\boldsymbol{\theta}_0,\ldots,\boldsymbol{\theta}_J)
          =\sum_{t=0}^{J-1}
            L(\boldsymbol{\theta}_t,\boldsymbol{\theta}_{t+1}).
          \label{eq:maxwell-discrete-action}
        \end{equation}

\section{Curvature action for Spidron-type origami with polygonal boundaries}
    \label{sec:spidron-boundary}
    
\begin{figure}[tb]
          \centering
          \includegraphics[width=\textwidth]
            {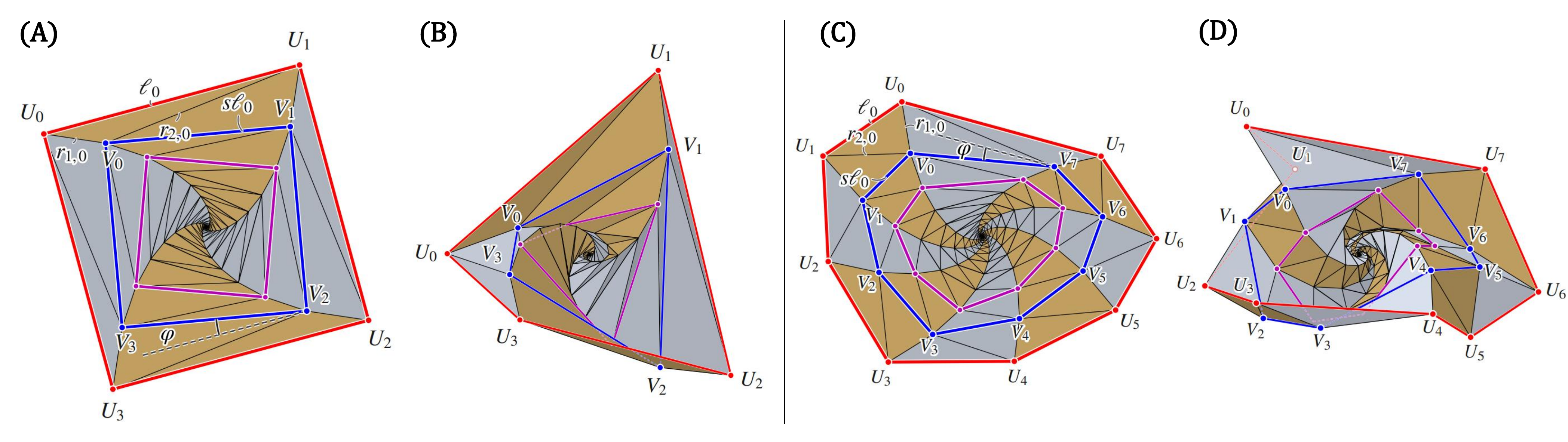}
            \caption{Spidron-type origami with quadrilateral and octagonal
            boundaries. (A) Crease pattern of the square-based Spidron studied
            by Imada et al. (B) Its folded configuration.
            (C) Crease pattern of an octagonal example with unequal side lengths.
            (D) Its folded configuration. Red, blue, and purple mark three
            successive boundaries, from the outside inward. In panels (A) and (C),
            the corresponding sides $U_0U_1$ and $V_0V_1$ have lengths $\ell_0$
            and $s\ell_0$, respectively, and $\varphi$ denotes the rotation angle
            between corresponding oriented sides.}
            \label{fig:spidron}
        \end{figure}

        In the preceding sections, boundary states propagate along a tube.
        We now show that the same curvature action applies to origami in which
        successive boundaries are nested in the crease pattern and boundary
        states propagate inward. The Spidron was introduced by
        Erd\'ely~\cite{Erdely2004Spidron}. Imada, Hull, Ku, and
        Tachi~\cite{Imada2026Spidron} studied a generalization whose crease
        pattern consists of nested squares. Each square is obtained from the
        preceding one by a rotation and a contraction about their common
        center, and the region between consecutive squares is triangulated.
        When the pattern is folded, the inner boundary of each ring serves as
        the outer boundary of the next. The boundary state thus propagates
        inward from one spatial quadrilateral to the next, as shown in
        Fig.~\ref{fig:spidron}(A), (B).
        
        Imada et al.~\cite{Imada2026Spidron} describe the shape of the outer
        quadrilateral by the exterior dihedral angles $\rho^1$ and $\rho^2$
        at its two diagonals. Following Section~\ref{subsec:maxwell_boundary_coordinates}, we choose a single
        diagonal, say $U_0U_2$, and take its half-length and the signed
        exterior dihedral angle at that same diagonal as our boundary
        variables. We then extend the construction from squares to simple polygons,
        as illustrated by the octagonal example in
        Fig.~\ref{fig:spidron}(C), (D), and derive the corresponding
        curvature action.

    \subsection{Polygonal rings and inward propagation}
        \label{subsec:spidron-ring-geometry}
        
        Fix a right-handed Cartesian coordinate system in three-dimensional
        Euclidean space. Place the crease pattern in the $xy$-plane and take
        its center of rotation and contraction, $O$, as the origin.
        Let $n\in\mathbb{Z}_{\geq 4}$, and let the outer boundary be a simple
        $n$-gon with vertices $U_0,\ldots,U_{n-1}$ in cyclic order.
        Throughout, vertex indices are taken in $\mathbb{Z}/n\mathbb{Z}$.
        As in Section~\ref{subsec:maxwell_boundary_coordinates}, boundary
        shapes are identified under translations and rotations, with the
        cyclic vertex labels fixed.
        
        Let $\mathbf{a}_i\in\mathbb{R}^2$ be the planar coordinates of $U_i$.
        Fix a contraction ratio $s\in(0,1)$ and a rotation angle
        $\varphi\in[-\pi,\pi)$, and let $R_\varphi\in\mathrm{SO}(2)$ denote
        rotation through $\varphi$. The planar coordinates
        $\mathbf{b}_i\in\mathbb{R}^2$ of the corresponding inner vertex $V_i$
        are given by
        \begin{equation*}
          \mathbf{b}_i=sR_\varphi\mathbf{a}_i.
        \end{equation*}
        We assume that the inner polygon lies strictly inside the outer polygon
        and that the segments $U_iV_i$ and $U_{i+1}V_i$ divide the intervening
        region into $2n$ nondegenerate triangles, as in
        Fig.~\ref{fig:spidron}(A), (C). We take this region as the reference
        ring. For $t\in\mathbb{Z}_{\geq 0}$, the $t$th ring of the crease
        pattern is obtained by rotating the reference ring through $t\varphi$
        about $O$ and scaling it by $s^t$. Since $0<s<1$, the planar rings
        converge to $O$ as $t\to\infty$.
        
        The crease pattern determines the edge lengths
        $\ell_i,r_{1,i},r_{2,i}\in\mathbb{R}_{>0}$ of the reference ring by
        \begin{equation*}
          \ell_i=\|\mathbf{a}_{i+1}-\mathbf{a}_i\|,\qquad
          r_{1,i}=\|\mathbf{b}_i-\mathbf{a}_i\|,\qquad
          r_{2,i}=\|\mathbf{b}_i-\mathbf{a}_{i+1}\|.
        \end{equation*}
        The outer edge $U_iU_{i+1}$ has length $\ell_i$, and the corresponding
        inner edge $V_iV_{i+1}$ has length $s\ell_i$. The edges $U_iV_i$ and
        $U_{i+1}V_i$ joining the two boundaries have lengths $r_{1,i}$ and
        $r_{2,i}$, respectively. All these lengths remain fixed as the ring
        folds.
        
        In a folded configuration, let $\mathbf{u}_i,\mathbf{v}_i\in\mathbb{R}^3$
        denote the position vectors of $U_i,V_i$, respectively, in the same
        coordinate system (Fig.~\ref{fig:spidron}(B), (D)). In the unfolded
        configuration, these vectors are obtained by appending a zero third
        coordinate to $\mathbf{a}_i$ and $\mathbf{b}_i$. Prescribe the position
        $\mathbf{u}_i$ of each outer vertex $U_i$ subject to
        $\|\mathbf{u}_{i+1}-\mathbf{u}_i\|=\ell_i$. We seek the position
        $\mathbf{v}_i$ of each inner vertex $V_i$ such that every triangle
        retains its prescribed edge lengths. Thus, for every $i$,
        \begin{equation}
          \|\mathbf{v}_i-\mathbf{u}_i\|=r_{1,i},\qquad
          \|\mathbf{v}_i-\mathbf{u}_{i+1}\|=r_{2,i},\qquad
          \|\mathbf{v}_{i+1}-\mathbf{v}_i\|=s\ell_i.
          \label{eq:spidron-folded-ring}
        \end{equation}
        
        The first two constraints in Eq.~\eqref{eq:spidron-folded-ring} concern
        a single triangle $[U_i,U_{i+1},V_i]$. With its edge $U_iU_{i+1}$
        fixed, the triangle is free only to rotate about that edge. The vertex
        $V_i$ therefore lies on a circle, on which we choose a local angular
        coordinate $\gamma_i\in\mathbb{R}$. The last constraint requires
        consecutive vertices $V_i,V_{i+1}$ to be joined by an inner edge of
        length $s\ell_i$. Once the outer vertex positions are prescribed, we
        solve these $n$ inner edge constraints simultaneously for the angles
        $\gamma_0,\ldots,\gamma_{n-1}$.
        
        The same outer boundary can admit more than one folded ring.
        For a quadrilateral ring, Imada et al.~\cite{Imada2026Spidron} give
        examples in which the same exterior dihedral angles $\rho^1,\rho^2$
        specifying the outer boundary admit distinct pairs of triangle
        inclinations $\phi^1,\phi^2$ compatible with the prescribed inner
        edge lengths. When several solutions exist, we select the one
        corresponding to the folding under consideration.
        
        Write the inner edge constraints as
        $G_i=\|\mathbf{v}_{i+1}-\mathbf{v}_i\|^2-s^2\ell_i^2=0$.

        The Jacobian $(\partial G_i/\partial\gamma_j)_{0\leq i,j<n}$
        vanishes at the planar crease pattern. We work locally near a selected
        folded configuration at which this Jacobian is nonsingular. The
        implicit function theorem then gives each $\gamma_i$ as a unique
        smooth function of the outer vertex positions, thereby determining
        the inner boundary shape.
        This defines a local propagation map from one boundary shape to the
        next. The inner boundary obtained at each step serves as the outer
        boundary of the next ring.
    
     \subsection{Oriented caps and the curvature action of one ring}
        \label{subsec:spidron-ring-action}
        
       \begin{figure}[tb]
          \centering
          \includegraphics[width=\textwidth]
            {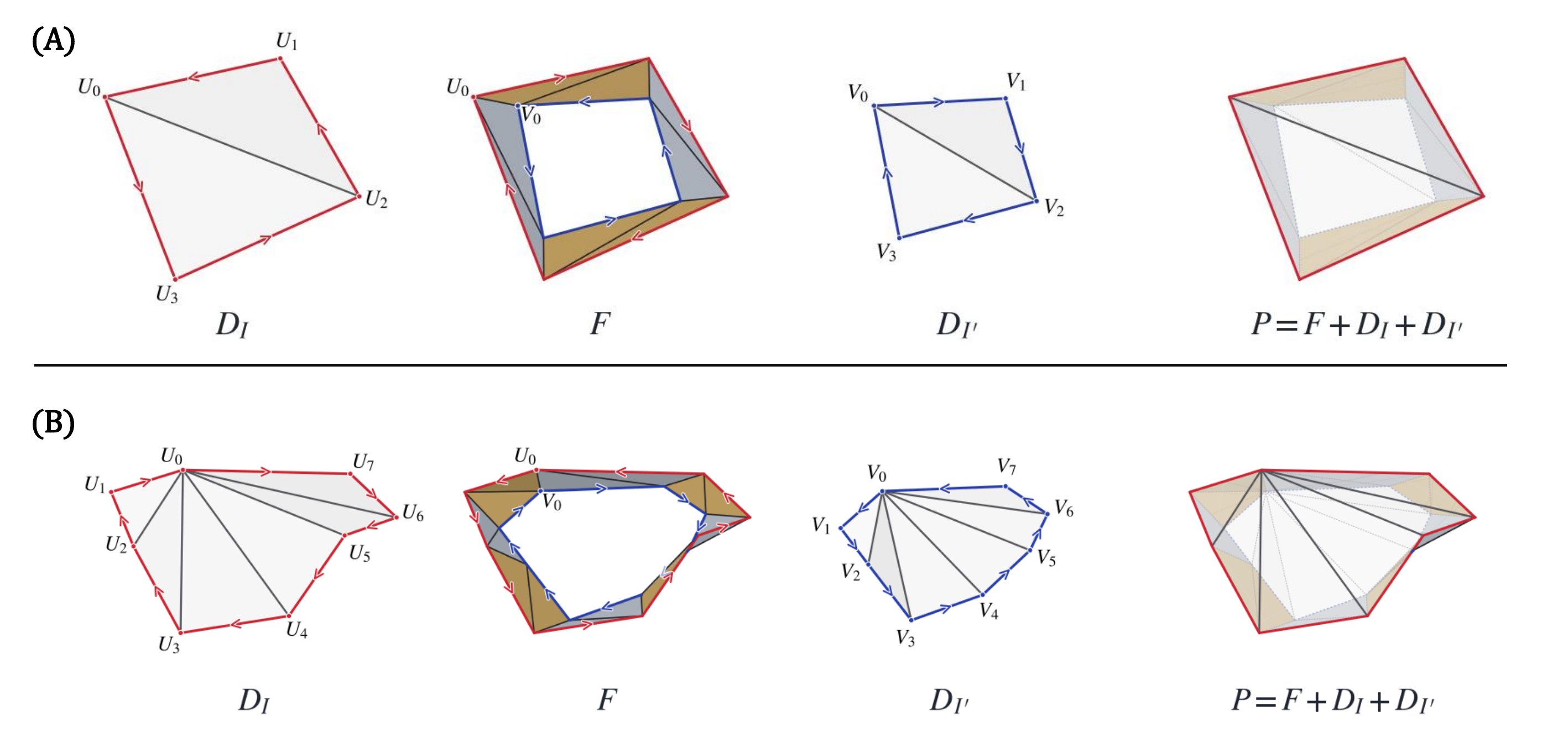}
            \caption{Construction of a closed oriented polyhedral surface from a
            single Spidron ring. (A) A ring with quadrilateral boundaries.
            (B) An octagonal example of the polygonal generalization introduced
            in this work. From left to right, each row shows the outer cap $D_I$,
            the folded ring $F$, the inner cap $D_{I^{\prime}}$, and their
            oriented sum $P=F+D_I+D_{I^{\prime}}$. Red and blue mark the outer
            and inner boundaries, respectively. Diagonals from $U_0$ and $V_0$
            triangulate the respective caps. The arrows show that the ring and
            each cap induce opposite orientations along their common boundary.
            Thus $\partial P=0$, and the Schl\"afli formula applies even when
            $P$ self-intersects.}
            \label{fig:spidron-caps}
        \end{figure}
        
        To construct a closed oriented polyhedral surface from a single ring,
        we introduce a cap for each of its two boundary polygons. Choose a
        triangulation of the outer polygon in the crease pattern. For each
        triangle of this triangulation, join the corresponding three vertices
        on the folded outer boundary and, likewise, on the folded inner
        boundary. The resulting triangles span the respective caps even when
        the boundary vertices are not coplanar. Fig.~\ref{fig:spidron-caps}(A)
        shows a quadrilateral ring, and Fig.~\ref{fig:spidron-caps}(B) shows
        an octagonal example of the polygonal generalization introduced here.
        In both examples, the outer and inner caps are triangulated by
        diagonals from $U_0$ and $V_0$, respectively.
        We restrict attention to configurations for which all triangles
        of the chosen input and output caps are nondegenerate.
        
        Orient both boundary polygons in order of increasing vertex index and
        write $C_U=\sum_{i=0}^{n-1}[U_i,U_{i+1}]$ and
        $C_V=\sum_{i=0}^{n-1}[V_i,V_{i+1}]$. As in Fig.~\ref{fig:closed-polyhedron}, denote the input
        and output caps by $D_I$ and $D_{I^{\prime}}$, respectively, and orient
        them so that $\partial D_I=-C_U$ and $\partial D_{I^{\prime}}=C_V$.
        The oriented sum $F$ of the physical faces of the ring and the
        resulting closed surface $P$ are given by
        \begin{equation*}
        \begin{aligned}
          F &= \sum_{i=0}^{n-1}
               \bigl([U_i,U_{i+1},V_i]+[U_{i+1},V_{i+1},V_i]\bigr),
          &\qquad \partial F &= C_U-C_V,\\
          P &= F+D_I+D_{I^{\prime}},
          &\partial P &= 0.
        \end{aligned}
        \end{equation*}
        The arrows in Fig.~\ref{fig:spidron-caps} show the opposite
        orientations induced by the ring and each cap along their common
        boundary. Intersections between faces leave the prescribed edge
        incidences unchanged, so the Schl\"afli formula of Proposition~\ref{prop:euclidean-schlafli}
        also applies when $P$ self-intersects.
        
        Following Section~\ref{subsec:maxwell_boundary_coordinates}, index the $n-3$ diagonals of each cap by
        $j=1,\ldots,n-3$. For the $j$th input diagonal, let
        $I_j\in\mathbb{R}_{>0}$ be its half-length and
        $\theta_j\in\mathbb{R}/(2\pi\mathbb{Z})$ the signed exterior dihedral
        angle between its two incident triangles. Denote the half-length and
        angle at the corresponding output diagonal by
        $I'_j\in\mathbb{R}_{>0}$ and
        $\theta'_j\in\mathbb{R}/(2\pi\mathbb{Z})$, respectively. For the
        variation below, choose locally smooth real-valued lifts of the
        angles, consistently with our angular convention, and collect
        the boundary variables into column vectors
        \[
        \begin{aligned}
          \mathbf I = (I_1,\ldots,I_{n-3})^{\mathsf T}, \  
          \boldsymbol\theta = (\theta_1,\ldots,\theta_{n-3})^{\mathsf T}, \
          \mathbf I' = (I'_1,\ldots,I'_{n-3})^{\mathsf T}, \
          \boldsymbol\theta' = (\theta'_1,\ldots,\theta'_{n-3})^{\mathsf T}.
        \end{aligned}
        \]
        
        To measure the boundary angles $\theta_j$ and $\theta'_j$, we orient
        each cap so that its boundary follows the increasing vertex indices.
        These orientations agree with those of $-D_I$ and $D_{I^{\prime}}$,
        respectively. Thus, on the oriented surface $P$, the exterior dihedral
        angle differentials at the input and output diagonals are
        $-\dd\theta_j$ and $\dd\theta'_j$, respectively.
        
        Let $E_{\mathrm{fix}}$ be the set of edges of the physical ring $F$,
        including the edges of both boundary polygons. For each
        $e\in E_{\mathrm{fix}}$, let $\ell_e$ be its length and
        $\Phi_e(\boldsymbol\theta,\mathbf I)$ its signed exterior
        dihedral angle on $P$, evaluated on the selected local solution
        with the reference crease pattern fixed. In the following geometric
        description, we suppress the arguments. At an interior edge of the
        ring, $\Phi_e$ is measured between two physical faces. At a boundary
        edge, it is measured between the physical face and the incident cap
        triangle. The curvature action $S_M$ of one ring and its differential
        are
        \begin{equation}
        \begin{aligned}
          S_M(\boldsymbol\theta,\mathbf I)
            =\frac12\sum_{e\in E_{\mathrm{fix}}}\ell_e
                \Phi_e(\boldsymbol\theta,\mathbf I), \qquad
          \dd S_M(\boldsymbol\theta,\mathbf I)
            =\mathbf I^{\mathsf T}\dd\boldsymbol\theta
              -(\mathbf I')^{\mathsf T}\dd\boldsymbol\theta'.
        \end{aligned}
          \label{eq:spidron-action-variation}
        \end{equation}
        
        Indeed, the lengths of the edges in $E_{\mathrm{fix}}$ remain fixed
        during folding, so their contribution to the sum in the Schl\"afli formula is
        $2\,\dd S_M$. The remaining edges are the cap diagonals. Their
        contributions are $-2\sum_{j=1}^{n-3}I_j\,\dd\theta_j$ at the input
        and $2\sum_{j=1}^{n-3}I'_j\,\dd\theta'_j$ at the output. By
        Proposition~\ref{prop:euclidean-schlafli}, the sum of these three contributions vanishes,
        which gives Eq.~\eqref{eq:spidron-action-variation}. In the next subsection,
        we normalize the scale of the output boundary and use this identity
        to construct the discrete action for a sequence of rings.

    \subsection{Scale normalization and the discrete variational principle of Spidron}
    \label{subsec:spidron-normalization}
    
    As in Subsection~\ref{sec:conformal}, we rescale the output boundary
    by $1/s$ to restore the reference side lengths $\ell_i$. For a
    polygonal boundary, this divides every diagonal half-length $I'_j$
    by $s$ and preserves the signed exterior dihedral angles $\theta'_j$.
    In the normalized coordinates, we use the same regular folding branch of the reference ring at every step.
    The normalized propagation map $\widehat M$ is thus defined locally
    on the shape space of polygons with side lengths $\ell_i$, with
    \[
      \widehat M(\boldsymbol\theta_t,\mathbf I_t)
      =(\boldsymbol\theta_{t+1},\mathbf I_{t+1})
      =(\boldsymbol\theta',\mathbf I'/s).
    \]
    
    Substituting $\mathbf I'=s\mathbf I_{t+1}$ into
    Eq.~\eqref{eq:spidron-action-variation} and taking the exterior derivative
    gives
    \begin{equation}
      \widehat M^*\omega=s^{-1}\omega,
      \qquad
      \omega=\sum_{j=1}^{n-3}\dd\theta_j\wedge \dd I_j.
      \label{eq:spidron-conformal-factor}
    \end{equation}
    Hence $\widehat M$ is locally conformally symplectic with factor
    $s^{-1}$ on the $(2n-6)$-dimensional boundary shape space, with one
    coordinate pair for each of the $n-3$ diagonals.
    
    Under the twist condition
    $\det(\partial\boldsymbol\theta_{t+1}/\partial\mathbf I_t)\ne 0$,
    the consecutive boundary angles $\boldsymbol\theta_t$ and
    $\boldsymbol\theta_{t+1}$ serve as local coordinates on the graph of
    $\widehat M$. Writing the curvature action of the reference ring in
    these coordinates as $L(\boldsymbol\theta_t,\boldsymbol\theta_{t+1})$,
    we obtain from Eq.~\eqref{eq:spidron-action-variation} the generating
    equations
    \begin{equation}
      \mathbf I_t
      =D_1L(\boldsymbol\theta_t,\boldsymbol\theta_{t+1}),
      \qquad
      \mathbf I_{t+1}
      =-s^{-1}D_2L(\boldsymbol\theta_t,\boldsymbol\theta_{t+1}).
      \label{eq:spidron-generating-equations}
    \end{equation}

    Since the curvature action is linear in the edge lengths at fixed angles,
    ring $t$ contributes $s^t$ times its reference-scale action. The discrete
    curvature action for $J\geq 1$ rings is therefore
    \begin{equation}
      \mathcal A(\boldsymbol\theta_0,\ldots,\boldsymbol\theta_J)
      =\sum_{t=0}^{J-1}s^t
        L(\boldsymbol\theta_t,\boldsymbol\theta_{t+1}).
      \label{eq:spidron-discrete-action}
    \end{equation}
    With $\boldsymbol\theta_0$ and $\boldsymbol\theta_J$ fixed, the
    stationarity condition $\delta\mathcal A=0$ is equivalent to the discrete
    Euler--Lagrange equations
    \begin{equation}
      D_2L(\boldsymbol\theta_{t-1},\boldsymbol\theta_t)
      +sD_1L(\boldsymbol\theta_t,\boldsymbol\theta_{t+1})
      =\mathbf 0,
      \qquad 1\leq t\leq J-1.
      \label{eq:spidron-discrete-euler-lagrange}
    \end{equation}
    By Eq.~\eqref{eq:spidron-generating-equations}, the boundary angle
    sequences of orbits of $\widehat M$ are precisely the stationary
    sequences of $\mathcal A$ in this coordinate neighborhood.
    
    At boundary $t$, the actual side lengths are $s^t\ell_i$ and the
    diagonal lengths are $2s^t I_{j,t}$. Together with the boundary angles,
    these lengths determine the boundary shape as in
    Subsection~\ref{subsec:maxwell_boundary_coordinates}. Given the placement
    of the initial boundary, the constraints Eq.~\eqref{eq:spidron-folded-ring},
    with lengths scaled by $s^t$ in ring $t$, determine the placement of
    each successive boundary on the selected regular solution.

\section{Numerical construction of folded origami configurations}
    \label{sec:numerical-examples}
    
    In this section, we numerically construct Maxwell tubes and
    Spidron-type origami with complex global shapes or distinctive
    symmetries. We characterize the resulting configurations as
    stationary solutions of the curvature action and illustrate how
    the variational structure can be used to design folded origami
    with prescribed geometric features.
    
    \subsection{Numerical method}
    \label{subsec:numerical-method}
    
    For a configuration of $J$ rings, the variables of the variational
    problem are the boundary angles
    $(\boldsymbol\theta_0,\ldots,\boldsymbol\theta_J)$.
    For the Maxwell tube, we use the discrete action in
    Eq.~\eqref{eq:maxwell-discrete-action}, which sums the curvature
    actions of the individual rings. For Spidron, we use
    Eq.~\eqref{eq:spidron-discrete-action}, in which the reference-scale
    curvature action of ring $t$ is weighted by $s^t$.
    Where the propagation map is smooth and the twist condition holds,
    the boundary-angle sequence of an orbit of $M$ for the Maxwell tube,
    or of $\widehat M$ for Spidron, is stationary for the corresponding
    action under variations that fix the endpoint angles. Conversely,
    a boundary-angle sequence satisfying the corresponding discrete
    Euler--Lagrange equations locally determines a propagation orbit.
    
    We first specify the number of rings $J$ and the crease-pattern
    parameters. For the Maxwell tube, these include the number of
    circumferential cells $N$ and the prescribed edge lengths
    $l_{\mathrm{L},n}$ and $l_{\mathrm{R},n}$, as described in
    Subsection~\ref{subsec:maxwell_propagation}. For Spidron, we use
    the number of outer boundary vertices $n$, their planar coordinates
    $\{\mathbf a_i\}$, the scale factor $s$, and the rotation angle
    $\varphi$ in the crease pattern, as described in
    Subsection~\ref{subsec:spidron-ring-geometry}. To specify the
    intended folded shape, we impose conditions on features such as
    the periodicity of the boundary states, the spatial placement
    or rotational symmetry of the boundaries, and the endpoint angles.
    
    We use the intended shape to construct an initial guess. Starting
    from this guess, we solve for a configuration that matches the
    prescribed edge lengths, joins adjacent rings, and satisfies the
    imposed design conditions. The following subsections summarize the geometric conditions imposed in the numerical examples.
    From the computed configuration, we obtain the boundary coordinates
    $(\boldsymbol\theta_t,\mathbf I_t)$. At each propagation step, we
    check the twist condition by computing the smallest singular value
    of the twist matrix. We check stationarity by evaluating the
    residuals of the discrete Euler--Lagrange equations, using the
    vector form of Eq.~\eqref{eq:discrete-euler-lagrange} for the
    Maxwell tube and Eq.~\eqref{eq:spidron-discrete-euler-lagrange}
    for Spidron.
    Writing $L_t=L(\boldsymbol\theta_t,\boldsymbol\theta_{t+1})$, we report
    \begin{equation}
      \varepsilon_{\mathrm{EL}}=
      \max_{t\in\mathcal T}\|D_2L_{t-1}+\lambda D_1L_t\|_\infty,
      \qquad
      \sigma_*=
      \min_{0\le t<J}\sigma_{\min}\!\left(
        \left.\frac{\partial\boldsymbol\theta_{t+1}}
                   {\partial\mathbf I_t}\right|_{\boldsymbol\theta_t}
      \right).
      \label{eq:numerical-el-residual}
    \end{equation}
    Here $\lambda=1$ for Maxwell tubes and $\lambda=s$ for Spidrons.
    The index set is $\mathcal T=\{1,\ldots,J-1\}$ for a finite segment
    and $\mathcal T=\mathbb Z/J\mathbb Z$ for a periodic orbit.
    
    The following examples satisfy the geometric constraints, the twist
    condition, and the stationarity conditions to numerical accuracy.
   The Supplementary Material provides the formulas for the curves used to construct initial configurations, the computational settings, the parameters of the computed solutions, and the details of these checks.

     \subsection{Folded configurations with prescribed geometric features}
        \label{subsec:prescribed-configurations}
        
\begin{figure}[p]
          \centering
          \includegraphics[
              width=\linewidth,
              height=0.80\textheight,
              keepaspectratio
            ]
            {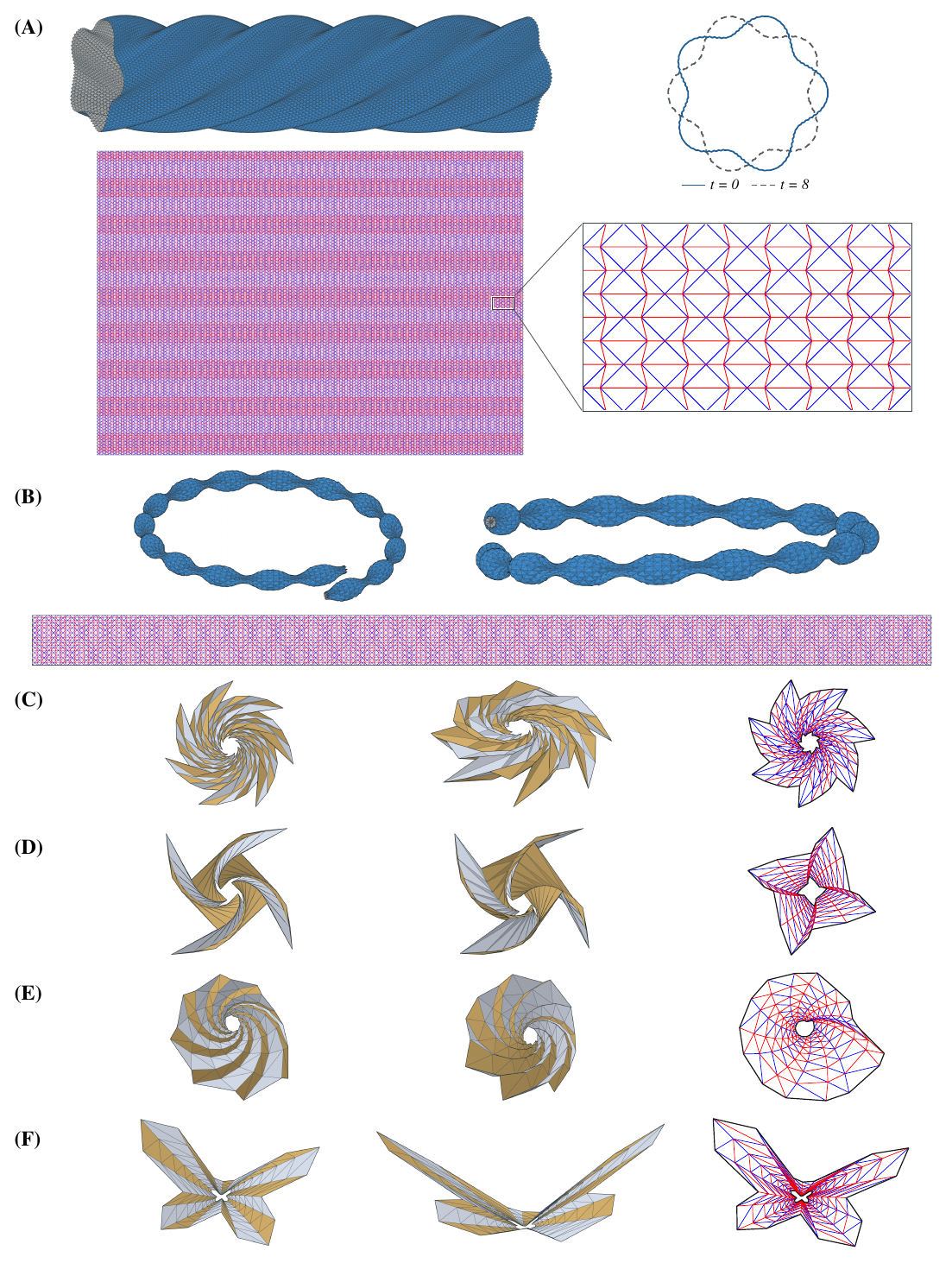}
            \caption{Folded configurations obtained numerically as stationary solutions
            of the discrete curvature action.
            (A) A Maxwell origami tube on the twist-mode branch with a prescribed
            five-lobed, star-shaped boundary (\(N=100\), \(J=80\)).
            The solid and dashed curves show the axial projections of the boundaries
            at \(t=0\) and \(t=8\), respectively. Below the tube are its complete
            crease pattern and an enlargement of the boxed region.
            (B) A coiling Maxwell origami tube with alternating bulges and necks
            (\(N=10\), \(J=144\)). Both views and the crease pattern show the first
            96 rings.
            (C)--(F) Spidron-type origami in the shapes of a seven-petalled flower,
            a four-bladed pinwheel, a spiral shell, and a butterfly, respectively.
            Each is shown from two directions, with its crease pattern on the right.
            Red and blue lines denote mountain and valley folds, respectively.
            No self-intersections were detected in any of the six folded surfaces.}
            \label{fig:prescribed-configurations}
        \end{figure}

        For the Maxwell origami tube in Fig.~\ref{fig:prescribed-configurations}(A),
        we prescribed a five-lobed boundary profile and obtained a stationary
        solution on the twist-mode branch. The boundary retains its star shape
        and rotates by approximately \(4.60^\circ\) per ring.
        For (B), we arranged the initial sequence of boundaries along a coil
        to guide the search. With the crease pattern and the boundary angles
        at both ends fixed, we obtained a tube that coils while alternately
        expanding and narrowing.
        
        For the Spidron-type configurations in (C)--(E), we sought successive
        folded boundaries that are similar in shape, with scale factor \(s\).
        These configurations correspond to fixed points of the normalized
        propagation map \(\widehat M\).
        For the butterfly shape in (F), we started from a shallow fold and
        obtained a more deeply folded stationary configuration with the same
        crease pattern and mountain and valley assignments.
    
    \subsection{Periodic configurations of Maxwell origami tubes}
        \label{subsec:periodic-maxwell-tubes}
        The propagation map \(M\) describes the shape of each boundary up to
        rigid motion. Along a periodic orbit of \(M\), the boundary therefore
        returns to its initial shape. For the tube to close in space, we also
        require each vertex of the final boundary \(C_J\) to coincide with the
        corresponding vertex of the initial boundary \(C_0\). To initialize the
        search, we place successive rings along a closed space curve chosen
        to have the desired knot type or projected shape.

\begin{figure}[p]
          \centering
          \includegraphics[width=\linewidth]
            {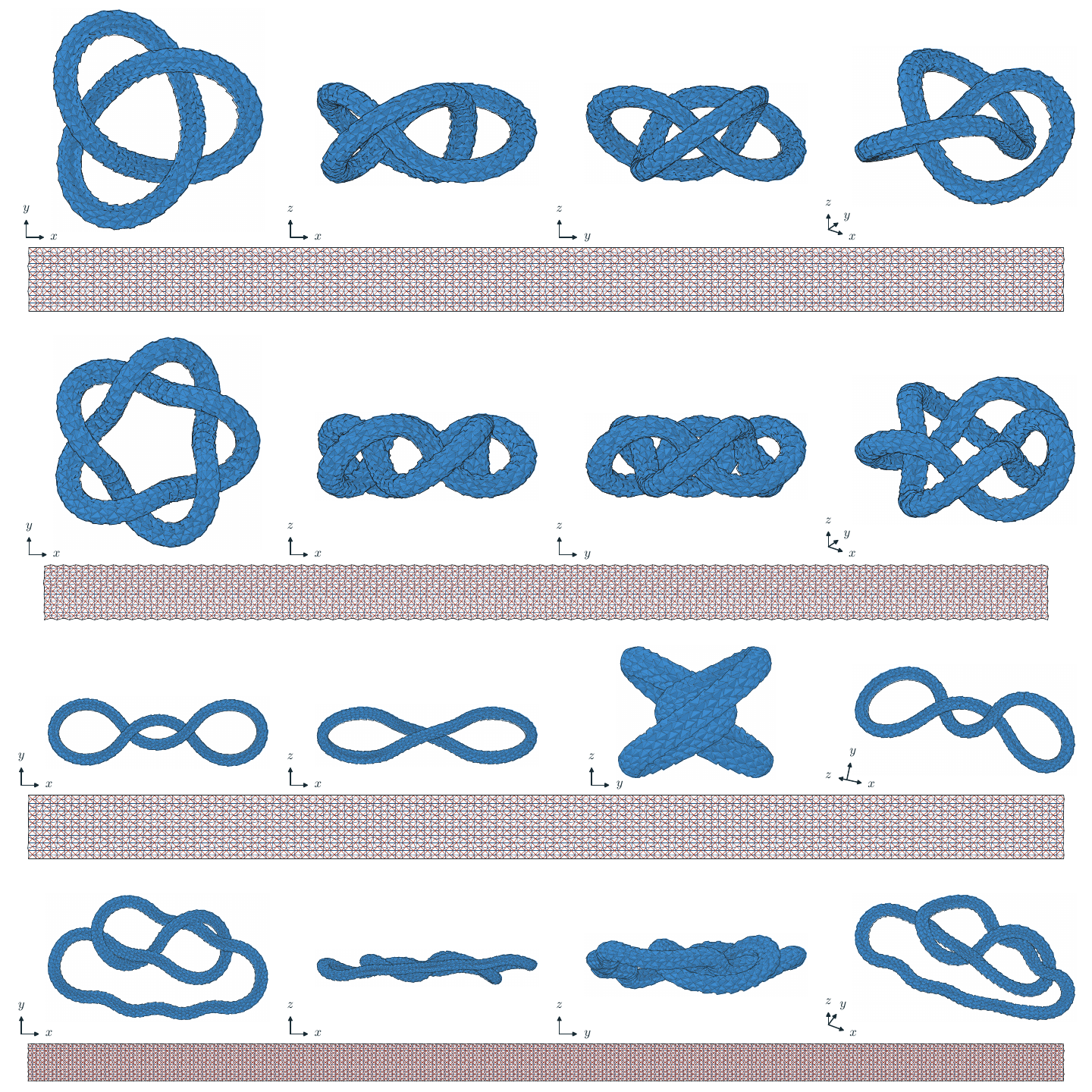}
            \caption{Periodic configurations of Maxwell origami tubes computed as
                stationary points of the discrete curvature action.
                From top to bottom, the centerlines form the trefoil knot \(3_1\),
                the cinquefoil knot \(5_1\), the unknot \(0_1\), and the figure-eight
                knot \(4_1\). The \(xy\)-projection of the third centerline resembles
                Gerver's super-eight orbit. Each tube has \(N=8\) circumferential cells.
                The tubes close in space after \(J=72\), \(80\), \(72\), and \(128\)
                iterations of the ring-to-ring propagation map, respectively.
                Each row shows, from left to right, projections onto the \(xy\)-,
                \(xz\)-, and \(yz\)-planes and an oblique view. The crease pattern
                below each row shows mountain folds in red and valley folds in blue.
                No self-intersections were detected in any of the four folded surfaces.}
              \label{fig:maxwell-periodic-tubes}
        \end{figure}

        For the first three configurations in Fig.~\ref{fig:maxwell-periodic-tubes},
        let \(b\in\mathbb{Z}_{\geq 1}\) denote the number of rings in one segment
        and \(k\in\mathbb{Z}_{\geq 2}\) the number of segments, so that \(J=bk\).
        Taking the origin on the rotation axis, we choose a rotation
        \(R\in\mathrm{SO}(3)\) of order \(k\) and impose, vertex by vertex,
        \[
          C_{t+b}=R(C_t), \qquad t\in\mathbb{Z}/J\mathbb{Z}.
        \]
        Thus each segment is rotated into the next, and the tube closes after
        \(J\) iterations. For the trefoil, cinquefoil, and super-eight
        configurations, we use \((b,k)=(24,3)\), \((16,5)\), and \((36,2)\),
        respectively. For the figure-eight knot, we impose spatial closure
        directly after \(J=128\) iterations without prescribing rotational symmetry.
        
        We evaluate the discrete Euler--Lagrange residual at every ring,
        including the junction between the last and first rings. From top to bottom in Fig.~\ref{fig:maxwell-periodic-tubes}, the maximum
        residuals \(\varepsilon_{\mathrm{EL}}\) are
        \(8.28\times10^{-13}\), \(4.62\times10^{-12}\),
        \(2.05\times10^{-13}\), and \(1.68\times10^{-10}\), respectively.

        \clearpage
        
        \clearpage

\section{Conclusion}
\label{sec:conclusion}

We constructed a discrete curvature action for the propagation of origami
boundary shapes directly from the folded geometry. When each ring is
closed with oriented auxiliary faces, the paper creases and the boundary
edges of the ring all have fixed lengths. The sum of the signed discrete
mean curvature contributions along these edges gives the one-ring action.
When adjacent rings are joined, their contributions along each shared
boundary edge combine to give the contribution associated with the signed
exterior dihedral angle between the two origami faces meeting at that edge.
Thus, the action for the full sequence of rings is the total signed
discrete mean curvature along the fixed-length edges of the two end
boundaries and all creases between them. Applying the Schl\"afli formula
to each ring yields a common variational identity that expresses the
variation of the one-ring action solely in terms of the input and output
boundary variables. Under a local twist condition, the one-ring action
serves as a generating function for the propagation map. The sequence
of boundary angles that determines the entire folded shape is therefore
a stationary point of the full curvature action when the endpoint angles
are held fixed. When successive rings differ in scale by a factor $s$,
this factor determines both the conformal factor $s^{-1}$ of the
normalized propagation map and the weight $s^t$ of the reference-scale
curvature action for ring $t$. Within this framework, we numerically
constructed knotted Maxwell origami tubes and polygonal Spidrons with
petal-like and pinwheel-like forms. For each example, we verified
geometric compatibility and the local twist condition at every step
and confirmed that the boundary-angle sequence satisfies the discrete
Euler--Lagrange equations to numerical accuracy. These checks confirm
that the computed shapes correspond to stationary sequences of the
curvature action and to orbits of their respective propagation maps.
These results show that the same variational principle can guide the
exploration of diverse folded shapes beyond tubular origami.

Some of the tubular origami configurations studied previously resemble constant mean curvature (CMC) surfaces when viewed at a scale larger than that of individual creases. The variational principle developed here requires the first variation of the curvature sum to vanish when the endpoint angles are held fixed. This stationarity condition allows the curvature contributions to vary across the folded surface. In fact, the knotted Maxwell origami tubes and the petal- and pinwheel-shaped polygonal Spidrons constructed here can all be described by stationary sequences under the same variational principle. How much freedom is there to vary the curvature distribution while preserving the prescribed edge lengths and compatibility between adjacent rings? For example, can we realize a prescribed overall shape by concentrating curvature in some regions while keeping the mean curvature nearly uniform in others? As the crease pattern is refined, which smooth surfaces can be approximated by such stationary configurations? The curvature action constructed here provides a way to investigate both the variety of surfaces that origami can realize and the limits imposed by its geometric constraints.
A related question is whether the curvature action extends to curved-crease origami through refinement of a polyhedral approximation.
As the number of boundary degrees of freedom increases, the curvature
action retains its geometric meaning as a sum of signed discrete mean
curvature contributions along fixed-length edges. This suggests a
possible extension to curved-crease origami. If the discrete action and
its variation converge under refinement of a polyhedral approximation
to a surface with curved creases, the limiting stationarity condition
could relate the local geometry of the creases to the shape of the
entire folded surface. A further open question is whether an analogous
curvature action can be constructed directly from the folded geometry
for origami dynamical systems beyond the symplectic and conformally
symplectic cases studied here, and whether its stationarity can
characterize propagation in those systems.

\bibliographystyle{unsrt}
\bibliography{references} 

\end{document}